\UseRawInputEncoding
\documentclass{article}
\usepackage{amsfonts}
\usepackage{amsmath}
\usepackage{amsthm}
\usepackage{color}
\usepackage[a4paper]{geometry}

\def\D{\Delta}

\def\i{\iota}
\def\l{\lambda}
\def\ra{\rightarrow}

\def\lr{\longrightarrow}

\def\o{\otimes}

\def\v{\varepsilon}
\def\vp{\varphi}

\theoremstyle{definition}
\newtheorem{definition}{Definition}[section]
\newtheorem{remark}[definition]{Remark}
\newtheorem{example}[definition]{Example}

\theoremstyle{plain}
\newtheorem{theorem}[definition]{Theorem}
\newtheorem{proposition}[definition]{Proposition}
\newtheorem{lemma}[definition]{Lemma}
\newtheorem{corollary}[definition]{Corollary}

\numberwithin{equation}{section}

\begin{document}
\title{\bf An algebraic framework for the Drinfeld double based on infinite groupoids }
\date{}
 \author{{\bf  Nan Zhou$^1$ and Shuanhong Wang$^{2}$\footnote {E-mail: nanzhou.math@zjsru.edu.cn(Nan Zhou), shuanhwang@seu.edu.cn(Shuanhong Wang)}}\\
 {\small 1. Department of Mathematics, Zhejiang Shuren University}\\
 {\small Hangzhou,Zhejiang 310015, China}\\
 {\small 2. School of Mathematics, Southeast University}\\
{\small Nanjing, Jiangsu 210096,  China}\\
}
 \maketitle

\begin{abstract}
The Drinfeld double associated to the weak multiplier Hopf ($*$-) algebra pairing $\left\langle  A, B\right\rangle$ is constructed. We show that the Drinfeld double is again a weak multiplier Hopf ($*$-) algebra. If $A$ and $B$ are algebraic quantum groupoids, then so does the double. We also prove  the correspondence between  modules over the Drinfeld double and Yetter-Drinfeld modules. Finally,  we prove that the double is a quasitriangular weak multiplier Hopf algebra.
\end{abstract}



\section{Introduction }

Weak multiplier Hopf algebra was introduced  by Van Daele and Wang in series papers \cite{VW1,VW2,VW3,VW4}. For any infinite groupoid, the groupoid algebra and its dual algebra are the basic examples of weak multiplier Hopf algebras.  And weak Hopf algebras (\cite{BNS,BS}) and multiplier Hopf algebras (\cite{VD,VD1}) are special cases of  weak multiplier Hopf algebras.  Weak multiplier Hopf algebra with  integrals is called  {\itshape algebraic quantum groupoid}. 

The aim of this paper is to give the Drinfeld double construction for algebraic quantum groupoid. Recall that the Drinfeld double  for finite dimensional weak Hopf algebras and multiplier Hopf algebras has been constructed in \cite{B1,NTV,ZhuW}  and \cite{De1,DeVD,DeVDW,DrVD,Zhang}, respectively. But the construction has not been considered yet for infinite weak (multiplier) Hopf algebras so far. From the double construction for multiplier Hopf algebras, we find that the pairing of two weak multiplier Hopf algebras  is a natural setting for the construction.    

The paper is organized as follows.

In Section 2  we recall the main definitions and properties on weak multiplier Hopf algebras.

In Section 3 we introduce a notation of   weak multiplier Hopf algebra pairing $\left\langle  A, B\right\rangle$. When there is a pairing $\left\langle  A, B\right\rangle$, we have  unital actions denoted by $\rhd $ and $\lhd$. We study a important map $R: A\o B \rightarrow M(A\o B)$.  By the Sweedler notation, $R$ is defined as
$$R(a\o b)=\sum b_{(1)}\rhd a\o b_{(2)}=\sum a_{(1)}\o b\lhd a_{(2)}.$$
Our  double will be constructed on the range of $R$.   

In Section 4 we give the main Drinfeld double construction.  Given a weak multiplier Hopf algebra pairing $\left\langle A, B \right\rangle$, we introduce the notation of Drinfeld double for weak multiplier Hopf algebras.  We will show that there exists a  weak multiplier Hopf algebra structure on the Drinfeld double algebra $D=A\o_D B$ in Theorem \ref{main}.  If $A$ and $B$ are algebraic quantum groupoids, then we get  a new algebraic quantum groupoid, see  Theorem \ref{Th1}.  In the case that $\left\langle A, B \right\rangle$ is a weak multiplier Hopf $*$-algebra pairing, we prove that $D$ is again a weak multiplier Hopf $*$-algebra. We also show that the left-left Yetter-Drinfeld module can be regarded as a left module over the Drinfeld double. This is a well know result in Hopf algebra theory.

In Section 5  we introduce the notion of quasitriangular weak multiplier Hopf algebra. We put the known results on the antipode of a quasitriangular Hopf algebra into the framework of weak multiplier Hopf algebras, see Proposition \ref{S^4}.  Finally, in Theorem \ref{QT-main} we show that the Drinfeld double $D$ is a quasitriangular weak multiplier Hopf algebra through the {\it canonical element}  in $M(A\o
 B)$. 
 \newline

{\bf Notation and basic references}
In this paper, we work with the field $\mathbb{C}$. $id$ means the identity map  and  $\tau$ means the flip map $A\o A \lr  A\o A, a\o b \mapsto b\o a$. The range of a map $f$ is denoted by $Range(f)$, the kernel of $f$ is denoted by $\ker(f)$. The composition of two maps $f$ and $g$ will be denoted by $f\circ g$, sometimes we write as $fg$.

Let $A$ be an algebra, $\mu_A$ means the multiplication in $A$. We always assume that the product $\mu_A$ is  {\itshape non-degenerate} as a bilinear form. $A^{op}$ means the opposite algebra. $A'$ is the space of all linear functionals on $A$.  $A$ is called {\itshape idempotent} if $A^2=A$. 

 The left multipliers of $A$ are maps $\rho_l : A \rightarrow A$ satisfying $\rho_l(ab)=\rho_l(a)b$ for all $a, b\in A$. Similarly elements of right multipliers are maps $\rho_r: A\rightarrow A$ satisfying $\rho_r(ab)=a\rho_r(b)$  for all $a, b\in A$.
The {\itshape multiplier algebra}  $M(A)$ of $A$ is defined as the set of pairs $(\rho_l, \rho_r)$ such that $\rho_r(a)b=a\rho_l(b)$ for all $a, b\in A$.    $M(A)$ can be  characterized as the largest algebra containing A as a dense ideal. Obviously, $M(A)$ is a unital algebra and if $A$ has an identity, then $M(A)=A$. 

Let $R$ be a ring. Given $a\in R$, if there exists $b\in R$ such that 
	$$aba=a\text{\quad and\quad}bab=b,$$
then $b$ is called the {\itshape generalized inverse} of $a$.  The notation  will be used frequently in the theory of weak multiplier Hopf algebras. Note that for the experts in generalized inverse theory, they prefer to call that $b$ is $\{1,2\}$-inverse of $a$, see \cite{AT}.

The theory of Hopf algebra can be found in \cite{Ma,M}. The use of the Sweedler notation for (weak) multiplier Hopf algebras can be found in \cite{VD2}.  For an  overview on the theory of weak multiplier  bialgebras we refer to \cite{B2,B3,BGL}.  

\section{Preliminaries}

In this section  we recall the definition of weak multiplier Hopf algebras and some basic properties from \cite{VW2,VW3,VW4,ZhouW}.

A {\it coproduct} on an  algebra $A$ is a homomorphism $\D: A\lr M(A\o A)$ such that
\begin{itemize}
 \item $\D(a)(1\o b)$ and $(a\o 1)\D(b)$ are in $A\o A$ for all $a, b\in A$,
 \item $\D$ is coassociative in the sense that
$$
(c\o 1\o 1)(\Delta\o \i)(\Delta(a)(1\o b))=(\i\o\Delta)((c\o 1)\Delta(a))(1\o 1\o b)
$$
for all $a, b, c\in A$.
\end{itemize}

The coproduct $\D$ is called \textit{full} if the smallest subspaces $V$ and $W$ of $A$ satisfying
$$\D(A)(1\o A)\subseteq V\o A \quad \quad \text{and} \quad \quad (A\o 1)\D(A)\subseteq A\o W$$
are equal to $A$. The coproduct $\D $ is called \textit{regular} if $\D(a)(b\o 1)$ and $(1\o a)\D(b)$ are in $A\o A$ for all $a, b\in A$.

The \textit{canonical maps} $T_1$ and $T_2$ from $A\o A$ to $A\o A$ are defined by
$$
T_1(a\o b)=\D(a)(1\o b) \quad \text{and} \quad T_2(a\o b)= (a\o 1)\D(b).
$$
If $\D$ is regular, we can define $T_3$ and $T_4$ on $A\o A$ by
$$
T_3(a\o b)=(1\o b)\D(a) \quad \text{and} \quad T_4(a\o b)=\D(b)(a\o 1).
$$

A linear map $\v: A\lr \mathbb{C}$ is called a {\it counit}
on an algebra  $A$ with a coproduct $\D $ if
$$
(\v\o id)(\D(a)(1\o b))=ab \quad \text{and} \quad (id \o \v)((a\o 1)\D(b))=ab
$$
for all $a, b\in A$.

Let $A$ be an algebra with a coproduct $\D$, define
 $\D ^{cop}: A\lr M(A\o A)$
via
$$
\D ^{cop}(a)(b \o c):= \tau(\D (a)(c \o b))\, \mbox { and}\,  (b \o c)\D ^{cop}(a):= \tau((c \o b)\D (a)).
$$
$A^{cop}$ is the algebra $A$ with the coproduct $\D^{cop}$.  

The leg numbering notation is used for elements and for maps. For example, $ \D _{13}: A \lr  M(A \o A \o A)$ is defined by
$$
\D _{13}(a)(b \o c \o d):= (id  \o \tau)(\D (a)(b \o d) \o c)
$$
 and
$$
(b \o c \o d)\D _{13}(a):= (id \o \tau)((b \o d)\D (a) \o c).
$$

Now we  recall the main definition from  \cite{VW2}.

\begin{definition}\label{wmha}
A weak multiplier Hopf algebra is a pair $(A,\D)$ of a non-degenerate idempotent algebra with a full coproduct and a counit satisfying the following conditions.

(i) There exists an idempotent $E\in M(A\otimes A)$ giving the ranges of the canonical maps:
$$
T_1(A\otimes A)=E(A\otimes A)\quad {and} \quad T_2(A\otimes A)=(A\otimes A)E.
$$

(ii) The element $E$ satisfies
$$
(id \otimes \D)E=(E\otimes 1)(1\otimes E)=(1\otimes E)(E\otimes 1),
$$

(iii) The kernels of the canonical maps are of the form
$$
Ker(T_1)=(1-G_1)(A\otimes A),$$
$$
 Ker(T_2)=(1-G_2)(A\otimes A),
$$
where $G_1$ and $G_2$ are linear maps from $A\otimes A$ to itself, given as
$$
(G_1\otimes id)(\D_{13}(a)(1\otimes b\otimes c))=\D_{13}(a)(1\otimes E)(1\otimes b\otimes c),
$$
$$
(id\otimes G_2)((a\otimes b\otimes 1)\D_{13}(c))=(a\otimes b\otimes 1)(E\otimes 1)\D_{13}(c),
$$
for all $a, b, c\in A$.
\end{definition}

A weak multiplier Hopf algebra is called \textit{regular} if the coproduct is regular and if
also $(A, \D^{cop})$ is a weak  multiplier Hopf algebra.
If $A$ is a $*$-algebra and $\Delta$ is a $*$-homomorphism, then we call $A$ a weak multiplier Hopf $*$-algebra.  Note that the regularity of  a weak  multiplier Hopf $*$-algebra is automatic. 

{\bf Notation.}
In the following,  weak multiplier Hopf algebra will be denoted by WMHA.

\begin{remark}
We will use the Sweedler notation in this paper. We will, e.g., write $\sum a_{(1)}\o a_{(2)}b$ for $\D(a)(1\o b)$. In the sequel  we just write $a_{(1)}\o a_{(2)}b$. The explanation of the Sweedler notation for these formulas can be found in \cite{VW3,ZhouW}. We are not going to explain it again. But in the following sections,  we always mention the necessary coverings  when the Sweedler notation was used.
\end{remark}

For any WMHA $A$, the {\it antipode} is defined as a unique map $S: A\ra M(A) $ which satisfies the following equalities
\begin{equation}\label{antipode}
\begin{split}
 & \sum S(a_{(1)})a_{(2)}S(a_{(3)})=S(a)\\
 & \sum a_{(1)}S(a_{(2)})a_{(3)}=a.
\end{split}
\end{equation}
  Now we can define two maps $P_1$ and $P_2$ from $A\o A$ to itself  by
\begin{align}\label{inverse}
P_1(a\o a')&=\sum a_{(1)} \o S(a_{(2)})a'\\
P_2(a\o a')&=\sum a S(a'_{(1)}) \o a'_{(2)}.
\end{align}
 They are called the generalized inverses of $T_1$ and $T_2$, respectively.

The following two important maps from $A$ to $M(A)$
$$\v_s: a\mapsto \sum S(a_{(1)})a_{(2)}, \quad \v_t: a\mapsto \sum a_{(1)}S(a_{(2)})$$
are called the {\itshape source} and {\itshape target maps} , respectively. Their ranges  $\v_t(A)$ and $\v_s(A)$ will be called the {\it target algebra} and {\it source algebra}, respectively. We can consider the multiplier algebras of $\v_t(A)$ and $\v_s(A)$. Denote $A_t=M(\v_t(A))$ and $A_s=M(\v_s(A))$, then
$$A_t=\{x\in M(A)|\D(x)=(x\o 1)E\} \quad  \quad A_s=\{y\in M(A)|\D(y)=E(1\o y)\}.$$
If $A$ is regular, we can define 
$$\v'_s(a)=\sum a_{(2)}S^{-1}(a_{(1)}) \text{\quad and \quad}\v'_t(a)=\sum S^{-1}(a_{(2)})a_{(2)}.$$

Now we  list some useful formulas.

\begin{lemma}\label{ASAT}
For any $x\in A_t, y\in A_s, a\in A$, we have

$(1)$ $\D(x)=(x\o 1)E=E(x\o 1) \quad \text{and} \quad \D(y)=E(1\o y)=(1\o y)E$,

$(2)$ $A_s$ and $A_t$ are commuting sub-algebras of $M(A)$,

$(3)$ $(1\o x)E=(S(x)\o 1)E \quad \text{and} \quad E(y\o 1)=E(1\o S(y))$,

$(4)$ $E(a\o 1)=\D(a_{(1)})(1\o S(a_{(2)}))\quad \text{and} \quad (1\o a)E=(S(a_{(1)})\o 1)\D(a_{(2)}).$

If $A$ is regular, we have

$(5)$ $E(1\o a)=\D(a_{(2)})(S^{-1}(a_{(1)})\o 1)\quad \text{and} \quad (a\o 1)E= (1\o S^{-1}(a_{(2)}))\D(a_{(1)})$.
\end{lemma}

The following theorem is  a characterization of weak multiplier Hopf algebras and it will be used in the proof of main double construction.

\begin{theorem}
Let $(A, \D)$ be a pair of  a non-degenerate idempotent algebra $A$ with a full coproduct $\D$ and such that there is a counit.

$(1)$ Assume that there is a linear map $S: A\ra M(A)$ such that the maps $P_1$ and $P_2$ defined above have range in $A\o A$, and such that for all $a\in A$,
$$a_{(1)}S(a_{(2)})a_{(3)}=a \quad \text{and} \quad S(a_{(1)})a_{(2)}S(a_{(3)})=S(a)$$
hold in $M(A)$.

$(2)$ Assume that there is an element $E$ in $M(A\o A)$ so that
$$T_1P_1(a\o b)=E(a\o b) \quad \text{and} \quad (T_2P_2)(a\o b)=(a\o b)E$$
for all $a, b\in A$ and that
$$(id\o \D)(E)=(\D\o id)(E)=(E\o 1)(1\o E)=(1\o E)(E\o 1).$$
Then $(A, \D)$ is a weak multiplier Hopf algebra and $S$ is the antipode.
\end{theorem}

Assume that $M$ is an $A$-module.  $M$ is called unital if $AM=M$. If $M$ is an algebra, then $M$ is called an $A$-module algebra if $M$ is a  unital $A$-module  and satisfies
 $$a(rr')=\sum(a_{(1)}\cdot r)(a_{(2)}\cdot r')$$
for all $a\in A, r, r'\in M$. Since the module is unital,  $a_{(1)}$ is well covered. And we have the following results.

\begin{proposition}
	Let M be a left A-module algebra. For any $r,r'\in M$ and $y\in \v_s(A)$, we have
	\begin{equation}\label{action}
	\mu_M(E\rhd (r\o r'))=rr'
\end{equation}
and
\begin{align}\label{abc}
(y\rhd r)r'=r(S(y)\rhd r'),
\end{align}
here $\mu_M$ is the multiplication in $M$.
\end{proposition}

 For more information about the action theory we refer to \cite{WZW, ZhouW}.

Let $(A, \D)$ be a regular WMHA. Now let us recall the notions of integral and dual algebra from \cite{VW4}. For any $a, b\in A$ and $f\in A'$, we can define  a multiplier $x\in M(A)$ by
$$xb=(f\o id)(\D(a)(1\o b))\quad \text{and} \quad bx=(f\o id)((1\o b)\D(a)).$$
And the multiplier $x$ is denoted by $(f\o id)\D(a)$. Similarly we can define multiplier $(id\o f)\D(a)$.

\begin{definition}
A non-zero linear functional $\vp: A\ra \mathbb{C}$ is called a {\it left integral} if $(id\o \vp)\D(a)\in A_t$ for all $a\in A$.  Similarly $\psi$ is called a {\it right integral} if $(\psi \o id)\D(a)\in A_s$ for all $a\in A$.
\end{definition}

Let $\int_L$ be the set of all left integrals on $A$ and $\int_R$ be the set of all right integrals on $A$.

\begin{definition}
	  We say that a WMHA $A$ has {\itshape a faithful set of integrals} if the following two conditions are satisfied. For any $x\in A$, we must have $x=0$ if $\vp(xa)=0$ for any $x\in \int_L$ and $a\in A$. Similarly, if  $\vp(ax)=0$ for any $x\in \int_L$ and $a\in A$, then $x=0$.
\end{definition}

\begin{definition}
	Let $(A, \D)$ be a regular WMHA. If $A$ has  a faithful set of integrals, then we call it an {\it algebraic quantum groupoid}.
\end{definition}

\begin{definition}
	Assume that $(A,\D)$ is an algebraic quantum groupoid. Then we define $\widehat{A}$ as the space of linear functionals on $A$ spanned by the elements of the form $\vp(\cdot a)$ where $\vp\in \int_L$ and $a\in A$.
\end{definition}

By the main theorem in \cite{VW4}, we know that $\widehat{A}$ is again an algebraic quantum groupoid with the dual structure. And we will call  it the dual algebraic quantum groupoid of $A$.

\section{Pairing}
\subsection{main definition and first properties }
Assume that $A$ and $B$ are two regular WMHAs. Let $\left\langle \cdot , \cdot \right\rangle  $ be a non-degenerate bilinear form from $A\times B$ to $\mathbb{C}$. For any  $a, b\in A$,  we have two maps
$$\left\langle a, \cdot \right\rangle  : B\ra \mathbb{C}, \quad b\mapsto \left\langle a, b\right\rangle  \quad\text{and}
\quad\left\langle \cdot, b \right\rangle  : A\ra \mathbb{C}, \quad a\mapsto \left\langle a, b\right\rangle.$$
 And next we can define the following four multipliers:
$$(\langle a,\cdot \rangle \o id)\D(b)\in M(B) \quad  \quad  (id\o \left\langle a, \cdot\right\rangle  )\D(b)\in M(B),$$
$$(\left\langle ,,b\right\rangle  \o id)\D(a)\in M(A)\quad \quad (id\o \left\langle \cdot, b\right\rangle  )\D(a)\in M(A).$$
We  use Sweedler notation here. For example, the element $ (\left\langle a, \cdot \right\rangle  \o id)\D(b)$ will be denoted by $\sum \left\langle a, b_{(1)}\right\rangle  b_{(2)}$, $b_{(1)}$ can be covered if we multiply an element belongs to $B$ from  right.

\begin{definition}\label{pairing}
Let $A$ and $B$ be two regular WMHAs. For any $a,a'\in A, b,b'\in B$,  the {\it pairing} between $A$ and $B$ is a non-degenerate bilinear form $\langle \cdot , \cdot \rangle$ from $A\times B$ to $\mathbb{C}$ satisfying the following:
\begin{itemize}
  \item[(1)]   $\sum a_{(1)}\left\langle a_{(2)}, b\right\rangle \in A $ and $\sum a_{(2)}\left\langle a_{(1)},b\right\rangle  \in A$,
  \item[(2)]  $\sum \left\langle a, b_{(1)}\right\rangle  b_{(2)} \in B$ and $\sum \left\langle a, b_{(2)}\right\rangle  b_{(1)} \in B$,
  \item[(3)] $\left\langle a, bb'\right\rangle  =\left\langle a_{(1)}, b\right\rangle  \left\langle a_{(2)}, b'\right\rangle  $,
  \item[(4)] $\left\langle aa', b,\right\rangle  =\left\langle a,b_{(1)}\right\rangle  \left\langle a', b_{(2)}\right\rangle  $
\end{itemize}
\end{definition}

Let $A$ and $B$ be two weak  multiplier Hopf $*$-algebras, if  $(A, B, \left\langle ,\right\rangle )$ is a weak multiplier Hopf algebra pairing with
\begin{equation}\label{*condition}
\left\langle a, b^*\right\rangle= \overline{\left\langle S(a)^*, b\right\rangle} \quad \text{and} \quad
\left\langle a^*, b\right\rangle= \overline{\left\langle a, S(b)^*\right\rangle},
\end{equation}
then we call $(A, B, \left\langle ,\right\rangle )$ a
{\itshape weak multiplier Hopf $*$-algebra pairing}.

Note that the items (1) and (2) of Definition \ref{pairing} give a meaning to the formulas in items (3) and (4).
For any pairing $\left\langle A, B\right\rangle$, we can define the following four  maps
\begin{eqnarray}
&&_{A}\rhd_B: A\o B\lr B: a\o b\mapsto \sum\left\langle a, b_{(2)}\right\rangle  b_{(1)},\label{map1}\\
&&_{B}\lhd_A : B\o A\lr B: b\o a\mapsto \sum\left\langle a, b_{(1)}\right\rangle  b_{(2)},\\
&& _{B}\rhd_A: B\o A\lr A: b\o a\mapsto \sum a_{(1)}\left\langle a_{(2)},b\right\rangle  ,\\
&&_{A}\lhd_B: A\o B\lr A: a\o b\mapsto \sum a_{(2)}\left\langle a_{(1)}, b\right\rangle  \label{map4}.
\end{eqnarray}

\begin{proposition}
	The above maps give four unital non-degenerate actions.
\end{proposition}
\begin{proof}
	 We only  consider the map $_{A}\rhd_B$. We claim that $B$ is spanned by the elements of  form $\{\sum\left\langle a, b_{(2)}\right\rangle  b_{(1)} \}$.  Suppose that this is not true, then there exists a non-zero linear functional $f$ on $B$ such that
	$$f(\sum \left\langle a, b_{(2)}\right\rangle  b_{(1)})=0\quad \quad \text{for all $a\in A, b\in B$}.$$
	Since the pairing is non-degenerate, we have $(f\o id)\D(b)=0$ in $M(A)$ for all $b$. By the fullness of coproduct, we get $f=0$, contradiction. So the action $_{A}\rhd_B$ is unital.
	Assume that $a\rhd b=\left\langle a, \sum b_{(2)}\right\rangle  b_{(1)}=0$ for all $a$. Then we have $\left\langle a,\sum b_{(2)}f(b_{(1)})\right\rangle =0$ for any $a\in A, f\in B'$. Again since  the coproduct  is full, we obtain $\left\langle a, b\right\rangle=0$ for all $b$ and it follows that $a=0$. Now we obtain a non-degenerate action.
\end{proof}

\begin{proposition}
Assume that $\left\langle A, B\right\rangle  $ is a WMHA pairing.  $(B, _{A}\rhd_B)$ is a left $A$-module algebra and $(B,_{B}\lhd_A)$ is a right $A$-module algebra. Analogously $(A,_{B}\rhd_A)$ is a left $B$-module algebra and $(A,_{A}\lhd_B)$ is a right $B$-module algebra.
\end{proposition}

For the proof we refer to the Proposition 3.13 in \cite{ZhouW}. These actions will be  denoted by $\rhd$ and $\lhd$.

\begin{proposition}\label{2}
Let $\left\langle A, B\right\rangle  $ be a WMHA pairing. For any $a, a'\in A, b, b'\in B$, we have
\begin{gather}
\left\langle aa', b\right\rangle  =\left\langle a, a'\rhd b\right\rangle  ,  \quad \quad \left\langle a'a, b\right\rangle  =\left\langle a,b\lhd a'\right\rangle    \\
\left\langle a, bb'\right\rangle  =\left\langle a\lhd b,  b'\right\rangle  ,  \quad \quad\left\langle a, b'b\right\rangle  =\left\langle b\rhd a, b'\right\rangle
\end{gather}
\end{proposition}

It is easy to prove.

\begin{proposition}
Let $\left\langle A , B\right\rangle  $ be a WMHA  pairing, then $(A, \rhd, \lhd)$ is a $B$-bimodule algebra and $(B, \rhd, \lhd)$ is a $A$-bimodule algebra.
\end{proposition}
\begin{proof}
For any $a\in A, b, b', b''\in B$, by the formulas in Proposition \ref{2} we have
$$\left\langle (b\rhd a)\lhd b', b''\right\rangle  =\left\langle (b\rhd a), b'b''\right\rangle  =\left\langle a, b'b''b\right\rangle  =\left\langle b\rhd (a\lhd b'), b''\right\rangle .$$
The rest of the argument is easy.
\end{proof}

Inspired by the Proposition 2.4  in \cite{VW4}, the WMHA pairing $(A, B, \left\langle \cdot, \cdot \right\rangle )$ can be extended in the following way.

\begin{proposition}
Let $(A, B, \left\langle \cdot, \cdot \right\rangle )$ be a WMHA pairing. Then there is a unique extension of the pairing of $A$ with $B$ to a pairing of $M(A)$ with $B$, satisfying
$$\left\langle m , b\lhd a\right\rangle  =\left\langle am, b\right\rangle    \quad \text{and}  \quad \left\langle m, a\rhd b\right\rangle  =\left\langle ma, b\right\rangle  $$
for any $a\in A, b\in B, m\in M(A)$.  Similarly, the pairing can be extended to $A$ with $M(B)$, satisfying
$$\left\langle a\lhd b , m\right\rangle  =\left\langle a, bm\right\rangle    \quad \text{and}  \quad \left\langle b\rhd a, m\right\rangle  =\left\langle a, mb\right\rangle  $$
for any $a\in A, b\in B, m\in M(B)$.
\end{proposition}

\begin{proof}
	Let us consider the first case.  Take any $m\in M(A), b\in B$. Since the action $\lhd $ and $\rhd $ are unital, we can write $b=\sum_i b_i\lhd a_i$. Now we will define
	$$\left\langle m, b\right\rangle  =\sum_i \left\langle a_im, b_i\right\rangle  .$$
	If $b=\sum_i b_i\lhd a_i=0$, take a local unit $e\in A$ such that $a_i=a_ie$, then
	$$\left\langle m, b\right\rangle  =\left\langle m, b_i\lhd a_i\right\rangle  =\left\langle m, b_i\lhd (a_ie)\right\rangle  =\left\langle em, b_i\lhd a_i\right\rangle  =0.$$
	It means that the definition is well defined.  For the second equation, we only need to  write $b=\sum_i b_i\rhd a_i$ and the proof is similar.  The uniqueness is obvious.
	
	For the second case, we can  try another way.  Take any $m\in M(B), a\in A$, we can write $a=\sum_i b_i\rhd a_i $. Now we will define
	$$\left\langle a, m\right\rangle  =\sum \left\langle a_i, mb_i\right\rangle  .$$
	If $a=\sum_i b_i\rhd a_i=0$, then for any $b\in B$ we have
	$$\left\langle (mb_i)\rhd a_i, b\right\rangle  =\left\langle a_i,b(mb_i)\right\rangle  =\left\langle b_i\rhd a_i, bm\right\rangle  =0.$$
	So $(mb_i)\rhd a_i=0$, apply counit $\v_A$ then we can get $\left\langle a_i, mb_i\right\rangle  =0$ .    Then we have showed that the definition is well defined. The proof of rest is easy.	
\end{proof}

\begin{remark}
(1) In the above proposition, we have used the result $\v(a\rhd b)=\left\langle a, b\right\rangle  $. If we  write $a=b_i\rhd a_i$, then
\begin{align*}
\v(a\rhd b) &= \v(\left\langle a,b_{(2)}\right\rangle   b_{(1)})\\
                 &= \v(\left\langle a_i, b_{(2)}b_{i}\right\rangle  b_{(1)}).
\end{align*}
Now $b_{(1)}$ is covered by $b_i$, so we can get
$$\v(a\rhd b)=\left\langle a_i, b_{(2)}b_{i}\right\rangle  \v(b_{(1)})=\left\langle a_i, bb_i\right\rangle  =\left\langle a, b\right\rangle .$$

(2) We can also get the extension of pairing from $(A\o A)\times (B\o B)$ to $(A\o A)\times M(B\o B)$ or $M(A\o A)\times (B\o B)$. But  the pairing can not be extended to $M(A)\times M(B)$ generally.
\end{remark}

For  any WMHA pairing $\left\langle  A,B\right\rangle $, now we can give an meaning to the expressions $\left\langle  a,1 \right\rangle  $ and $\left\langle 1,b\right\rangle  $ where $a\in A, b\in B$. See the following corollary.

\begin{corollary}
Let $(A, B, \left\langle \cdot, \cdot \right\rangle  )$ be a WMHA pairing, then we have
\begin{equation}
\left\langle a,1\right\rangle  =\v_A(a),\quad   \quad  \left\langle  1,b \right\rangle   =\v_B(b)
\end{equation}
for any $a\in A, b\in B$, here $1$ belongs to the corresponding multiplier algebra.
\end{corollary}

\begin{proof}
We will show the first formula, the second one can be proved symmetrically. Write $a=a_i\lhd b_i$, here $a, a_i\in A, b_i\in B$. According to the extension we have
$$\left\langle a,1\right\rangle  =\left\langle a_i, b_i1\right\rangle  =\left\langle a_i, b_i\right\rangle  =\v_A(a_i\lhd b_i)=\v_A(a).$$
\end{proof}

\begin{proposition}
Let $(A, B, \left\langle \cdot, \cdot \right\rangle  )$ be a WMHA pairing. For all $a,a'\in A, b, b'\in B$, we have
\begin{itemize}
	\item [(1)] $\left\langle a\o a', E^B\right\rangle  =\left\langle aa', 1\right\rangle  =\v_A(aa')$

    \item [(2)] $\left\langle E^A, b\o b'\right\rangle  =\left\langle 1, bb'\right\rangle  =\v_B(bb').$
\end{itemize}
\end{proposition}

\begin{proof}
	(1) For any $a, a'\in A, b,b'\in B$,  take $x, y\in B$ such that  $E^B(b\o b')=\D_B(x)(1\o y)$.
	\begin{align*}
	\left\langle b\rhd a \o b'\rhd a', E^B\right\rangle    &= \left\langle a\o a', E^B(b\o b')\right\rangle  \\
	          									&= \left\langle a\o a', \D_B(x)(1\o y)\right\rangle  \\
									            &= \left\langle aa'_{(1)}, x\right\rangle  \left\langle a'_{(2)}, y\right\rangle  \\
									            &=  \left\langle a(y\rhd a'), x\right\rangle  \\
									            &=   \left\langle (x_{(1)}\rhd a)(x_{(2)}y\rhd a'), 1\right\rangle  \\
									            &= \v_A((b\rhd a)(b'\rhd a'))
	\end{align*}
In the fifth equality we use the definition of module algebra. In the last equality we use the formula (\ref{action}). Since the action is unital, we finish the proof.

  (2) The proof is similar.
\end{proof}

Note that the formulas in  Proposition 3.9 are used to define the canonical idempotent $\widehat{E}$ on the dual space  $\widehat{A}$ (see Proposition 2.13 in \cite{VW4}).

  Let $\left\langle A,B\right\rangle  $ be a WMHA pairing, for any $a\in A, b\in B$, by  Lemma 1.4 in \cite{DrVD} we know that
$$(id \o \cdot \lhd a)\D_B(b), \quad (id \o a\rhd \cdot)\D_B(b)$$
$$(\cdot \lhd a\o id )\D_B(b), \quad (a\rhd \cdot \o id)\D_B(b)$$
are multipliers in $M(B\o B)$, and
$$(id \o \cdot \lhd b)\D_A(a), \quad (id \o b\rhd \cdot)\D_A(a)$$
$$(\cdot \lhd b\o id )\D_A(a), \quad (b\rhd \cdot \o id)\D_A(a)$$
are multipliers in $M(A\o A)$.
Now we will give some formulas which are related to these multipliers.

\begin{proposition}\label{aaa}
Let $\left\langle A, B\right\rangle  $ be a pairing, for any $a, b\in A$, we have
\begin{description}
  \item[(1)] $\D_B(a\rhd b)=(id\o a\rhd \cdot)\D_B(b)$
  \item[(2)] $\D_B(b\lhd a)=(\cdot\rhd a\o id)\D_B(b)$
  \item[(3)] $\D_A(a\lhd b)=(\cdot\lhd b\o id)\D_A(a)$
  \item[(4)] $\D_A(b\rhd a)=(id\o b\rhd \cdot)\D_A(a)$
\end{description}
\end{proposition}
\begin{proof}
Using the coassociativity of $\D$, it is easy to prove these formulas.
\end{proof}

\begin{remark}
	
Note that for any $a\in A, b\in B, m\in M(B)$,  the action $A\lhd B$ can be extended. So we can define $a\lhd m$. Since we have
$$a\lhd (bm)=\left\langle a_{(1)}, bm\right\rangle  a_{(2)}=\left\langle a_{(1)}\lhd b,m\right\rangle  a_{(2)}=\left\langle (a\lhd b)_{(1)}, m\right\rangle  (a\lhd b)_{(2)}.$$
In the last equality we use the item (3) of Proposition \ref{aaa}. Now $a\lhd m$ can be denoted as $\left\langle a_{(1)}, m\right\rangle  a_{(2)}$.  Moreover,  we can  give  meaning to the formula $$\left\langle a, mn\right\rangle  =\left\langle a_{(1)}, m\right\rangle  \left\langle a_{(2)}, n\right\rangle  ,$$ here $a\in A, m,n\in M(B)$.
\end{remark}

\begin{proposition}
Let $\left\langle A, B\right\rangle  $ be a WMHA pairing, for any $a, a'\in A, b, b'\in B$, we have
\begin{description}
  \item[(1)] $\left\langle T^A_1(a\o a'), b\o b'\right\rangle  =\left\langle a\o a',T^B_2(b\o b')\right\rangle  $
  \item[(2)] $\left\langle T^A_2(a\o a'), b\o b'\right\rangle  =\left\langle a\o a',T^B_1(b\o b')\right\rangle  $
\end{description}
\end{proposition}
\begin{proof}
(1) For the left side we have
\begin{align*}
\left\langle T^A_1(a\o a'), b\o b'\right\rangle   &= \left\langle \D(a)(1\o a'),b\o b'\right\rangle  \\
                          &= \left\langle a_{(1)},b\right\rangle  \left\langle a_{(2)}a', b'\right\rangle  \\
                          &\overset{(3.4)}{=} \left\langle a\lhd b, a'\rhd b'\right\rangle  .
\end{align*}
For the right side we have
\begin{align*}
\left\langle a\o a',T^B_2(b\o b')\right\rangle   &= \left\langle a\o a', (b\o 1)\D(b')\right\rangle  \\
                         &= \left\langle a\o a', bb'_{(1)}\o b'_{(2)}\right\rangle  \\
                         &\overset{(3.1)}{=} \left\langle a,b(a'\rhd b')\right\rangle  \\
                         &\overset{(3.6)}{=} \left\langle a\lhd b, a'\rhd b'\right\rangle  .
\end{align*}
(2) Similar.
\end{proof}

\begin{proposition}
Let $\left\langle A, B\right\rangle  $ be a WMHA pairing,  the following equations hold.
\begin{align}
   & \left\langle \v^A_t(a), b\right\rangle  =\left\langle a,\v^B_t(b)\right\rangle  \\
   & \left\langle \v^A_s(a), b\right\rangle  =\left\langle a,\v^B_s(b)\right\rangle  \\
   & \left\langle S_A(a),b\right\rangle  =\left\langle a,S_B(b)\right\rangle  \label{S(a)baS(b)}.
\end{align}
\end{proposition}
\begin{proof}
 Recall that
$$\v_t(a)=(\v\o id)(E(a\o 1)) \quad \text{and} \quad \v_s(a)=(id\o \v)((1\o a)E).$$
For the formula (3.9), we have
\begin{align*}
\left\langle \v^A_t(a), b\right\rangle   &= \left\langle  \v^A(E^A_1a)E^A_2, b \right\rangle   \\
                 &\overset{(3.7)}{=} \left\langle  E^A_1a, 1 \right\rangle   \left\langle  E^A_2, b \right\rangle  \\
                 &= \left\langle  E^A_1,E^B_1\right\rangle  \left\langle  a, E^B_2\right\rangle  \left\langle  E^A_2, b \right\rangle
\end{align*}
In the last equality we use $\D_B(1)=E_B=E^B_1\o E^B_2$. If we write $a=ea$, where $e$ is a local unit, then $E^A_1$ can be covered by the local unit $e$. Thus the formula in the last equality is meaningful.

On the other hand we have
\begin{align*}
\left\langle a, \v^B_t(b)\right\rangle   &= \left\langle a, \v^B(E^B_1b)E^B_2\right\rangle   \\
                 &\overset{(3.7)}{=} \left\langle  1,E^B_1b \right\rangle   \left\langle  a, E^B_2\right\rangle  \\
                 &= \left\langle  E^A_1,E^B_1\right\rangle  \left\langle  E^A_2, b\right\rangle  \left\langle  a, E^B_2\right\rangle
\end{align*}
 The last equality can be explained in similar way. So we can get $\left\langle \v^A_t(a), b\right\rangle  =\left\langle a, \v^B_t(b)\right\rangle .$

 The proof of (3.10) is similar.

 Now let us prove the final equation. Take any $a, a'\in A, b, b'\in B$, then for the left side we have
 \begin{align*}
 \left\langle S_A(a\lhd b'), a'\rhd b\right\rangle   &\overset{(\ref{antipode})}{=} \left\langle S_A((a\lhd b')_{(1)})(a\lhd b')_{(2)}S_A((a\lhd b')_{(3)}), a'\rhd b\right\rangle  \\
               &= \left\langle a_{(1)}, b'\right\rangle   \left\langle \v^A_s(a_{(2)}), b_{(1)}\right\rangle   \left\langle S_A(a_{(3)})a', b_{(2)}\right\rangle  \\
               &\overset{(3.9)}{=} \left\langle a_{(1)}, b'\right\rangle   \left\langle a_{(2)}, \v^B_s(b_{(1)})\right\rangle   \left\langle S_A(a_{(3)})a', b_{(2)}\right\rangle  \\
               &= \left\langle a_{(1)}, b'\right\rangle   \left\langle a_{(2)}, S_B(b_{(1)})\right\rangle   \left\langle a_{(3)}, b_{(2)}\right\rangle   \left\langle S_A(a_{(4)})a', b_{(3)}\right\rangle  \\
               &= \left\langle a_{(1)}, b'\right\rangle   \left\langle a_{(2)}, S_B(b_{(1)})\right\rangle   \left\langle \v^A_t(a_{(3)}), a'\rhd b_{(2)}\right\rangle  \\
               &\overset{(3.8)}{=} \left\langle a_{(1)}, b'\right\rangle   \left\langle a_{(2)}, S_B(b_{(1)})\v^B_t(a'\rhd b_{(2)})\right\rangle  \\
               &= \left\langle a_{(1)}, b'\right\rangle   \left\langle a_{(2)}, S_B(a'\rhd b)\right\rangle  \\
               &= \left\langle a\lhd b', S_B(a'\rhd b)\right\rangle
 \end{align*}
In the second equality, $a_{(1)}$ and $a_{(3)}$ are covered by $b'$ and $a'$ respectively, then the formula is meaningful. In the fourth equality we use the
\begin{equation}
\left\langle a,\v_s(b)\right\rangle  =\left\langle a_{(1)}, S(b_{(1)})\right\rangle  \left\langle a_{(2)}, b_{(2)}\right\rangle  \label{bbb}
\end{equation}
If we write $a$ as $b'\rhd a$, then
\begin{align*}
\left\langle b\rhd a, \v_s(b)\right\rangle  &=   \left\langle a, S(b_{(1)})b_{(2)}b'\right\rangle  \\
                              &=\left\langle a_{(1)}, S(b_{(1)})\right\rangle  \left\langle a_{(2)}, b_{(2)}b' \right\rangle \\
                              &=\left\langle a_{(1)}, S(b_{(1)}) \right\rangle \left\langle b'\rhd a_{(2)}, b_{(2)} \right\rangle\\
                               &= \left\langle (b'\rhd a)_{(1)}, S(b_{(1)}) \right\rangle  \left\langle (b'\rhd a)_{(2)}, b_{(2)} \right\rangle.
\end{align*}
So the formula \eqref{bbb} is meaningful. In the fifth and seventh equalities we use the result in Proposition \ref{aaa}. In the last equality we use the definition of action $\lhd$.

Since the two actions $_{B}\rhd_A$ and $_{B}\lhd_A$ are unital, we prove the formula (3.9).
\end{proof}

 As an application of formula (3.11), it is easy to obtain the following  results.
\begin{corollary}
Let $\left\langle A, B\right\rangle  $ be a WMHA pairing,  for any $a\in A, b\in B$ the following equations hold.
\begin{description}
  \item[(1)] $S_A(b\rhd a)=S_A(a)\lhd S_B^{-1}(b)$
  \item[(2)] $S_B(a\rhd b)=S_B(b)\lhd S^{-1}_A(a),$
\end{description}
\end{corollary}
\begin{proof}
    \item[(1)] Recall that for a regular WMHA $(A,\D ,S)$, we have
	$\D(S(a))=\tau (S\o S)\D(a)$
	for any $a\in A$.  Now
	\begin{align*}
	S_A(a)\lhd S_B^{-1}(b)=S(a_{(1)})\left\langle S(a_{(2)}), S^{-1}(b)  \right\rangle= S(a_{(1)})\left\langle a_{(2)}, b  \right\rangle=S_A(b\rhd a).
	\end{align*}
	\item[(2)] The proof is similar.
\end{proof}

Recall that $P_1$ and $P_2$ are the generalized inverses  of the canonical maps $T_1$ and $T_2$, respectively.

\begin{proposition}
Let $\left\langle A, B\right\rangle  $ be a WMHA pairing, for any $a, a'\in A, b, b'\in B$, we have
\begin{description}
  \item[(1)] $\left\langle P^A_1(a\o a'), b\o b'\right\rangle  =\left\langle a\o a',P^B_2(b\o b')\right\rangle  $
  \item[(2)] $\left\langle P^A_2(a\o a'), b\o b'\right\rangle  =\left\langle a\o a',P^B_1(b\o b')\right\rangle  $
\end{description}
\end{proposition}
\begin{proof}
(1)
\begin{align*}
\left\langle P^A_1(a\o a'), b\o b'\right\rangle    &\overset{(\ref{inverse})}{=} \left\langle (id\o S_A\o id)(\D_A(a)\o a'), b\o \D_B(b')\right\rangle  \\
                          &\overset{(3.10)}{=} \left\langle \D_A(a)\o a', b\o (S\o id)\D_B(b')\right\rangle  \\
                          &= \left\langle  a\o a', (b\o 1)(S_B\o id)\D(b')\right\rangle  \\
                          &\overset{(\ref{inverse})}{=} \left\langle a\o a', P^B_2(b\o b')\right\rangle
\end{align*}

(2) The proof is  similar.
\end{proof}

Since we have $E=T_1P_1,G_2=T_2P_2$, then we can get the following corollary.
\begin{corollary}
Let $\left\langle A, B\right\rangle $ be a WMHA pairing, for any $a, a'\in A, b, b'\in B$, we have
$$
\left\langle E^A(a\o a'), b\o b'\right\rangle =\left\langle a\o a', G^B_2(b\o b')\right\rangle
$$
and
$$
\left\langle G^B_2(a\o a'), b\o b'\right\rangle =\left\langle a\o a', E^B(b\o b')\right\rangle .
$$
\end{corollary}

\subsection{The map $R$}

Let $(A, B, \left\langle \cdot, \cdot\right\rangle  )$ be a WMHA pairing, define the map $R: A\o B \ra M(A\o B)$  as
\begin{equation}
\begin{split}
&R(a\o b)(a'\o b')=(id\o \left\langle \cdot, \cdot\right\rangle  \o id)(\D_A(a)(a'\o 1)\o \D_B(b)(1\o b'))\\
&(a'\o b')R(a\o b)=(id\o \left\langle \cdot, \cdot\right\rangle  \o id)((a'\o 1)\D_A(a)\o(1\o b')\D_B(b))
\end{split}
\end{equation}
here $a, a'\in A, b, b'\in B$. Similarly we can define the map
$$\widetilde{R}:=(id\o \left\langle \cdot, \cdot\right\rangle  \o id)(\D_B\o \D_A): B\o A\ra M(B\o A).$$

Since the coproduct can be extended to the multiplier algebra, the map $R$ can be defined from $M(A)\o M(B)$ to $M(A\o B)$.  
For any $x\in A\o B, a\in A, b\in B, \vp\in B', \psi \in A'$, the set
$$\{(id\o \vp)(R(x)(a\o b))\}$$
is called  {\it left leg} of $R$, the set
$$\{(\psi\o id)((a\o b)R(x))\}$$
is called  {\it right leg} of $R$.  Similarly,
$\{(id\o \psi)(\widetilde{R}(x)(b\o a))\}$
is called the {\it left leg} of $\widetilde{R}$ and
$\{(\vp\o id)((b\o a)\widetilde{R}(x))\}$
is called the {\it right leg} of $\widetilde{R}$.

\begin{definition}
	The map $R$ is called \it{full} if the left leg of $R$ is $A$ and the right leg of $R$  is $B$. Similarly, $\widetilde{R}$ is called \it{full} if its  left  and right legs are  $B$ and  $A$ respectively.
\end{definition}

\begin{proposition}\label{full}
	Let $(A, B, \left\langle \cdot, \cdot\right\rangle  )$ be a pairing, then   $R$ and $\widetilde{R}$ are full.
\end{proposition}
\begin{proof}
	For any $a, a' \in A, b, b' \in B, \vp\in B'$,  by the definition of $R$,
	\begin{align}\label{mapR}
	R(a\o b)(a'\o b')= (b_{(1)}\rhd a)a'\o b_{(2)}b'=a_{(1)}a'\o (b\lhd a_{(2)})b'.
	\end{align}
	 Since $\D$ is full and the action is unital,  we can get that $R$ is  full. The proof of $\widetilde{R}$ is similar.
\end{proof}

 Recall the smash product in \cite{ZhouW},  we use similar notation $E_B\rhd (A\o B)$  and $(A\o B)\lhd E_A$. The element will be denoted by $E^B_1\rhd a\o E^B_2b$ and $aE^A_1\o b\lhd E^A_2$, respectively. The expression is meaningful because the action is unital. From equation (\ref{mapR}), we get the following result about the range of $R$.

\begin{corollary}
	We have $Range(R)=E^B\rhd (A\o B)=(A\o B)\lhd E^A$.
\end{corollary}

	 In (multiplier) Hopf algebras case,   $R$ is a bijective map.  Now for weak multiplier Hopf algebras,  we will show that $\widetilde{R}$ and $R$ are  generalized inversed.

\begin{proposition}
Let $(A, B, \left\langle \cdot, \cdot\right\rangle  )$ be a  WMHA pairing, there exist   linear maps $R'$ and $\widetilde{R}'$ such that
$$R'RR'=R' \quad \text{and} \quad RR'R=R$$
$$\widetilde{R}'\widetilde{R}\widetilde{R}'=\widetilde{R}' \quad \text{and} \quad \widetilde{R}\widetilde{R}'\widetilde{R}=\widetilde{R}.$$
\end{proposition}
\begin{proof}
Define
$$R':= (id\o \left\langle \cdot, \cdot\right\rangle   (S^{-1}_A\o id)\o id)(\D_A\o \D_B),$$
it is from $A\o B$ to $M(A\o B)$. Now let us consider the range of $R'$. For any $a, a'\in A, b, b', b''\in B$, we have
$$
(a'\o 1)R'(b\rhd a\o b')(1\o b'')= (a'\o 1)(S^{-1}_B(b'_{(1)})b\rhd a\o b'_{(2)})(1\o b'').
$$
It means that
$$R'(_{B}\rhd_A\o id)=(_{B}\rhd_A \o id)(\tau\o id)(id\o (S^{-1}_B\o id)T^B_2(S_B\o id))(\tau\o id)$$
Denote $_{B}\rhd_A \o id$ by $f$ and $(\tau\o id)(id\o (S^{-1}_B\o id)T^B_2(S_B\o id))(\tau\o id)$ by $h$, then we have
$$R'f=fh.$$

On the other hand, we consider the range of $R$. From the proof of Proposition \ref{full}, we have
$$Rf=fg,$$
here
$g= (\tau\o id)(id\o \tau T^{B,op}_1 \tau)(\tau\o id).$
By a standard  computation, we can get
$$ghg=g \quad \text{and} \quad hgh=h.$$
Now we have
$$
R'RR'f=R'Rfh=R'fgh=fhgh=fh=R'f,
$$
$$
RR'Rf=RR'fg=Rfhg=fghg=fg=Rf.
$$
Since $f$ is surjective, it follows that $R'RR'=R'$ and $RR'R=R$.

$\widetilde{R}'$ is defined as
$$\widetilde{R}':= (id\o \left\langle \cdot ,\cdot \right\rangle   (S^{-1}_A\o id)\tau \o id)(\D_B\o \D_A)): B\o A\ra B\o A.$$
 The  proof is similar.
\end{proof}

Denote $N=RR', Q=R'R$, then $N$ and $Q$ are idempotent maps. Moreover,  we have
 \begin{align*}
 (a'\o 1)P(b\rhd a\o b')(1\o b'') &= (a'\o 1)(b'_{(2)}S^{-1}_B(b'_{(1)})b\rhd a\o b'_{(3)})(1\o b'')\\
                                             &= (a'\o 1)(E^B_1b\rhd a\o E^B_2b')(1\o b'').
 \end{align*}
 It means that
 $$  Range(R)=Range(N)=E^B\rhd (A\o B).$$
 Note that $RR'R=R$, then we have
 $$\ker(R)=Range (1-Q).$$

\begin{proposition}
	The generalized inverses of $R$ and $\widetilde{R}$ are unique.
\end{proposition}
\begin{proof}
	We have found that the $N$ and $Q$ are idempotent maps so that $N$ projects on the range $R$ and $1-Q$ projects on the kernel of $R$. By Lemma 2.1 of \cite{VW2},  the uniqueness can be deduced. The proof of $\widetilde{R}$ is completely similar. 
\end{proof}

\begin{proposition}
	$R$ and $\widetilde{R}$ satisfy, for all $a, a', b, b'$,
\begin{itemize}
	\item[$(1)$] $R(a\o bb')=b\rhd R(a\o b')$
	\item[$(2)$] $(\D_A\o id)R=(id\o R)(\D_A\o id)$\quad and \quad $(id\o \D_B)R=(R\o id)(\D\o id)$
	\item [$(1')$] $\widetilde{R}(b\o aa')=a\rhd \widetilde{R}(b\o a')$
	\item [$(2')$] $(\D_B\o id)\widetilde{R}=(id\o \widetilde{R})(\D_B\o id)$ \quad and \quad $(id\o \D_A)\widetilde{R}=(\widetilde{R}\o id)(id\o \D_A)$
\end{itemize}
Also the map $(id\o R)$ and $(\widetilde{R}\o id)$ commute with each other.
\end{proposition}
\begin{proof}
	(1) $R(a\o bb')=b_{(1)}\rhd (b'_{(1)}\rhd a)\o b_{(2)}b'_{(2)}=b\rhd R(a\o b').$
	
	(2) The two formulas follow from the coassociativity of the coproduct.
	
	The proof of $(1')$ and $(2')$ is easy. Finally,
	\begin{align*}
	(id\o R)(\widetilde{R}\o id)(b'\o a\o b)&=(1\o R)(a_{(1)}\rhd b'\o a_{(2)}\o b)\\
	&= a_{(1)}\rhd b'\o b_{(1)}\rhd a_{(2)}\o b_{(2)}\\
	&=(\widetilde{R}\o id)(id\o R)(b'\o a\o b) .
	\end{align*}
	In the second equality, $a_{(1)}$ can be covered if denote $b'$ as $a'\rhd b'$.
\end{proof}
\begin{remark}
As we see, the maps $R$ and $\widetilde{R}$ have  much in common with  canonical map $T_1$ and $T_2$. Their ranges are both determined by the idempotent element $E$ and their kernels have the similar expressions. They are both generalized inversed. They have the same behaviors with respect to coproduct. 
$R$ and $\widetilde{R}$ will play an important role in the next section. They also lead us to consider the twist tensor product of weak multiplier Hopf algebras, it will be discussed in another paper.
\end{remark}

Now we will give some basic examples of weak multiplier Hopf algebra pairing.

\begin{example}
Let $A$ and $B$ be usual Hopf algebras.  then the dual pairing (see Definition 1.4.3 in \cite{Ma}) between them is a weak multiplier Hopf algebra pairing. 
\end{example}

\begin{example}
If $A$ and $B$ are multiplier Hopf algebras,  then the multiplier Hopf algebra pairing $(A, B, \left\langle \cdot, \cdot\right\rangle  )$ (see Definition 2.8 in \cite{DrVD}) is  a  weak multiplier Hopf algebra pairing.  
\end{example}

\begin{example}
If $A$ and $B$ are weak Hopf algebras, then the weak Hopf dual pairing (see Definition 3.1 in \cite{ZhuW}) between $A$ and $B$ is a  weak multiplier Hopf algebra pairing. 
\end{example}

\begin{example}
Let $A$ be a  regular weak multiplier Hopf algebra with a faithful set of integrals,  then the dual $\widehat{A}$ is again a regular weak multiplier Hopf algebra with integrals (see Theorem 2.15 in \cite{VW4}). Define $\left\langle a, f\right\rangle  :=f(a)$, where $a\in A, f\in \widehat{A}$. Then $(A, \widehat{A}, \left\langle \cdot, \cdot\right\rangle  )$ is a non-degenerate weak multiplier Hopf algebra pairing. For the verification we refer to \cite{VW4} where the the pairing $\langle A,\widehat{A}\rangle  $ were investigated. From the Proposition 2.3 in \cite{VW4} we get four unital actions. For example, the right action of $\widehat{A}$ on $A$ is given as $a\lhd f=f(a_{(1)})a_{(2)}$ and this right action is unital. The Proposition 2.5 and 2.10 in \cite{VW4} give us $\left\langle a, 1\right\rangle  =\v_A(a)$ and $\left\langle 1,f\right\rangle  =\widehat{\v}(f)$.  The antipode $\widehat{S}$ of $\widehat{A}$ satisfies
$\langle a, \widehat{S}(b)\rangle  =\left\langle S(a),b\right\rangle  $ for all $a\in A, b\in \widehat{A}$.
\end{example}

Next we  consider the $*$-case.

\begin{proposition}
	Let $(A, B, \left\langle \cdot, \cdot\right\rangle )$ be a weak multiplier Hopf $*$-algebra pairing, we have
		  $$(a\rhd b)^*=S(a)^*\rhd b^* \quad\text{and}\quad (b\lhd a)^*=b^*\lhd S(a)^*,$$
		  $$(b\lhd a)^*=S(b)^*\lhd a^* \quad\text{and}\quad (a\rhd b)^*=a^*\rhd S(b)^*,$$
	for any $a\in A, b\in B$.	
\end{proposition}
\begin{proof}
	  For any $a\in A, b,b'\in B$, we compute
	  \begin{align*}
	  (a\rhd b)^*b' &= (\left\langle a, b_{(2)}\right\rangle b'^*b_{(1)})^*\\
	  &= \left\langle S^{-1}(a^*), b_{(2)}^*\right\rangle b_{(1)}^*b'\\
	  &= (S^{-1}(a^*)\rhd b^*)b'=(S(a)^*\rhd b^*)b',
	  \end{align*}
	  here $b'$ is used for the coverings. Since the product is non-degenerate, we obtain the first equation. All the other equations can be proved in a similar way.
\end{proof}

\begin{proposition}
	Let $\left\langle A, B \right\rangle$  be a weak multiplier Hopf $*$-algebra pairing, then we have
	$$R\circ(*\o *)=(*\o *)\circ R'  \quad R'\circ (*\o *)=(*\o *)\circ R$$
	$$R^{op,cop}\circ (*\o *)=(*\o *)\circ R^{op,cop,'} \quad
	R^{op,cop,'}\circ(*\o *)=(*\o *)\circ R^{op,cop}.$$
\end{proposition}
\begin{proof}
	For any $a, a'\in A, b, b'\in B$, then
	\begin{align*}
	(a'\o b')R(*\o *)(a\o b) &= a'a^*_{(1)}\left\langle a^*_{(2)}, b^*_{(1)}\right\rangle\o b'b^*_{(2)}\\
	&=a'a^*_{(1)}\overline{\left\langle a_{(2)}, S^{-1}(b_{(1)})\right\rangle}\o b'b^*_{(2)}\\
	&=a'a^*_{(1)}\o b^*_{(2)}\overline{\left\langle S^{-1}(a_{(2)}), b_{(1)}\right\rangle}\\
	&=(a'\o b')(*\o *)R'(a\o b).
	\end{align*}
	The proof of rest is similar.
\end{proof}

As we know that if $A$ is a regular weak multiplier Hopf $*$-algebra, then $A^{op, cop}$ is again a regular weak multiplier Hopf $*$-algebra with the same canonical idempotent and antipode. So we can get the following result.

\begin{proposition}
Let $(A, B, \left\langle \cdot, \cdot\right\rangle )$ be a weak multiplier Hopf $*$-algebra pairing, then\\ $(A^{op,cop}, B^{op,cop},\left\langle \cdot, \cdot\right\rangle )$ is again a weak multiplier Hopf $*$-algebra pairing.
\end{proposition}
\begin{proof}
	We have to check the four conditions in Definition \ref{pairing}. The first two conditions are obviously. We also have
	$$\left\langle a, b\cdot_{op}b'\right\rangle=\left\langle a, b'b\right\rangle=\left\langle a_{(1)}, b'\right\rangle \left\langle a_{(2)}, b\right\rangle=\left\langle a_{(2)}, b\right\rangle\left\langle a_{(1)}, b'\right\rangle,$$
	similarly we can get $\left\langle a \cdot_{op} a', b\right\rangle=\left\langle a, b_{(2)}\right\rangle\left\langle a', b_{(1)}\right\rangle$. Since the antipode of $A^{op,cop}$ is the same with $A$, the equation (\ref{*condition}) holds. Then we finish the proof.
\end{proof}

For the WMHA pairing $(A^{op,cop}, B^{op,cop},\left\langle \cdot, \cdot\right\rangle  )$, we can define
$$R^{op,cop}:=(id\o \left\langle \cdot, \cdot\right\rangle  \o id)(\D^{cop}_A\o \D^{cop}_B): A^{op,cop}\o B^{op,cop}\ra M(A^{op,cop}\o B^{op,cop}).$$
The map $\widetilde{R^{op,cop}}$ from $B^{op,cop}\o A^{op,cop}$ to $M(B^{op,cop}\o A^{op,cop})$ can be defined similarly.  Their generalized inverses will be denoted by $R^{op,cop,'}$ and $\widetilde{R^{op,cop,'}}$, respectively.

\begin{proposition}
Let $(A, B, \left\langle \cdot, \cdot\right\rangle  )$ be a WMHA pairing, then the following equalities hold.
\begin{description}
\item[(1)] $R^{op,cop} \circ (S^{\pm 1}\o S^{\mp 1})=(S^{\pm 1}\o S^{\mp 1}) \circ R$
\item[(2)] $R\circ R^{op,cop} = R^{op,cop}\circ R.$
\end{description}
\end{proposition}
\begin{proof}
(1) We start from the left side.
	\begin{align*}
	R^{op,cop} \circ (S^{\pm 1}\o S^{\mp 1})
	 &= (id\o \left\langle \cdot, \cdot\right\rangle  \o id)(\D^{cop}_A\o \D^{cop}_B)(S^{\pm 1}\o S^{\mp 1}) \\
	 &= (id\o \left\langle \cdot, \cdot\right\rangle  \o id)((S^{\pm 1}\o S^{\pm 1})\o (S^{\mp 1}\o S^{\mp 1}))(\D_A\o \D_B)\\
	 &= (S^{\pm}\o \left\langle \cdot, \cdot\right\rangle \o S^{\mp})(\D_A\o \D_B)\\
	 &=(S^{\pm 1}\o S^{\mp 1}) (id \o \left\langle \cdot, \cdot\right\rangle \o id)(\D_A\o \D_B)\\
	 &= (S^{\pm 1}\o S^{\mp 1}) \circ R
	\end{align*}
	In the second equality we use $\D^{cop}\circ S=(S\o S)\D$.
	
(2) For any $a, a', a''\in A, b\in B$, we have
\begin{align*}
R^{op,cop} R(a\o a'\rhd b\lhd a'') &= R^{op,cop}(a_{(1)}\o a'\rhd b\lhd (a''a_{(2)}))\\
                								  &=a_{(2)}\o (a_{(1)}a')\rhd b\lhd (a''a_{2})\\
                								  &= R(a_{(2)}\o (a_{(1)}a')\rhd b \lhd a'')\\
                								  &= RR^{op,cop}(a\o a'\rhd b\lhd a'')
\end{align*}
The elements $a'$ and $a''$ are used for the coverings.
\end{proof}

For any WMHA pairing $\left\langle A, B \right\rangle$, we give a characterization of elements in $M(B)$ in the end of  this section.

\begin{proposition}
	Let $\left\langle A, B \right\rangle$ be a WMHA paring. Given a linear functional $\omega\in A'$. Then there is  a multiplier $m\in M(B)$ such that $\left\langle a, m \right\rangle=\omega(a)$ if and only if
	$$(\omega\o id)\D(a)\in A \quad \text{and} \quad (id\o \omega)\D(a)\in A$$
	for all $a\in A$.
\end{proposition}
\begin{proof}
	(1) Assume that we have such multiplier $m$, then
	$$(\omega\o id)\D(a)=\sum\left\langle a_{(1)}, m \right\rangle a_{(2)}=a\lhd m,$$
	$m$ is covered if we replace $a$ by $\sum_i b_{i}\lhd a_{i}$ and $a\lhd m$ is in $A$. Similarly for the second expression.
	
	(2) Conversely, for any $b\in B$, let $x, y\in B$ such that
	$$\left\langle a, x\right\rangle= \omega(b\rhd a) \quad \text{and} \quad \left\langle a, y\right\rangle=\omega(b\rhd a).$$
	Define a multiplier $m$ by $mb=x$ and $bm=y$.  For  any $a\in A, b, b'\in B$, we have
	$$\left\langle a, b'(mb) \right\rangle=\left\langle a\lhd b', mb \right\rangle=\omega(b\rhd (a\lhd b'))$$
	$$\left\langle a, (b'm)b \right\rangle=\left\langle b\rhd a, b'm \right\rangle=\omega((b\rhd a)\lhd b')$$
	This means that $b'(mb)=(b'm)b$, so $m$ is a multiplier of $B$.
\end{proof}

\section{The Drinfeld double of algebraic quantum groupoid }

In this section  we give  the main Drinfeld double construction.

\subsection{The algebra structure of Drinfeld double}

Let $(A, B, \left\langle \cdot, \cdot\right\rangle  )$ be a non-degenerate weak multiplier Hopf algebra pairing. Define a twist map as
\begin{equation}
T:= R\circ R^{op,cop,'}\circ \tau: B\o A\ra A\o B.
\end{equation}
Define a product on $A\o B$ as
$$\cdot_D:= (\mu_A\o \mu_B)(id\o T\o id).$$

\begin{proposition}
 The multiplication $\cdot_D$ is associative.
\end{proposition}

\begin{proof}
In fact it is enough to show that the twist map $T$ satisfies the following two equations:
\begin{align}
T\circ(\mu_B\o id)=(id\o \mu_B)(T\o id)(id\o T)\\
T\circ (id\o \mu_A)=(\mu_A\o id)(id\o T)(T\o id).
\end{align}
 We will prove the first one, the second one can be derived in a similarly way. For any $a\in A, b, b', b'', b'''\in B$,
\begin{align*}
& T(m_B\o id)(b'\o b''\o a\lhd b)(1\o b''') \\
&= T(b'b''\o a\lhd b)(1\o b''')\\
&= RR^{op,cop,'}(a\lhd b\o b'b'')(1\o b''')\\
&= R(a\lhd s^{-1}((b'b'')_{(2)}S(b))\o (b'b'')_{(1)})(1\o b''')\\
&= (b'b'')_{(1)}\rhd a\lhd S^{-1}((b'b'')_{(3)}S(b))\o (b'b'')_{(2)}b'''
\end{align*}
On the other side, we have
$$T(b\o a\lhd b')(1\o b'')=b_{(1)}\rhd a\lhd S^{-1}(b_{(3)}S(b'))\o b_{(2)}b''$$
and
$$T(b\o b'\rhd a)(1\o b'')=b_{(1)}b'\rhd a\lhd S^{-1}(b_{(3)})\o b_{(2)}b''.$$
Then
\begin{align*}
& ((id\o m_B)(T\o id)(id\o T)(b'\o b''\o a\lhd b))(1\o b''') \\
&= (id\o m_B)(T\o id)(b'\o  b''_{(1)}\rhd a\lhd S^{-1}(b''_{(3)}S(b))\o b''_{(2)}b''')\\
&= (b'b'')_{(1)}\rhd a\lhd S^{-1}((b'b'')_{(3)}S(b))\o (b'b'')_{(2)}b'''
\end{align*}
\end{proof}

\begin{proposition}
The multiplication $\cdot_D$ is well-defined in $R(A\o B)$.
\end{proposition}
\begin{proof}
Recall that
$R(A\o B)=E^B\rhd (A\o B)=(A\o B)\lhd E^A$. 	
For any $a, a'\in A, b, b'\in B$, we have
\begin{align*}
(a\o b)\cdot_D(E^B_1\rhd a'\o E^B_2b') &= a(b_{(1)}\rhd (E^B_1\rhd a')\lhd S^{-1}(b_{(3)}))\o b_{(2)}E^B_2b'\\
                                &= a(b_{(1)}\rhd a'\lhd S^{-1}(b_{(3)}))\o b_{(2)}b'\\
                                &= (a\o b)\cdot_D(a'\o b')
\end{align*}
Note that the formulas above are meaningful if we write $a'$ as $a''\rhd a\lhd b''$ (since the actions are unital). In this way elements are well covered. In the second equality we use $\D(b)=\D(b)E^B$.

On the other hand, we have
\begin{align*}
(E^B_1\rhd a'\o E^B_2b')\cdot_D(a\o b) &= (E^B_1\rhd a')((E^B_2b')_{(1)}\rhd a \lhd S^{-1}((E^B_2b')_{(3)}))\o (E^B_2b')_{(2)}b\\
                                &= (E^B_1\rhd a')(E^B_2b'_{(1)}\rhd a \lhd S^{-1}(b'_{(3)}))\o b'_{(2)}b\\
                                &= a'(b'_{(1)}\rhd a \lhd S^{-1}(b'_{(3)}))\o b'_{(2)}b\\
                                &= (a'\o b')\cdot_D(a\o b)
\end{align*}
In the second equality we use
$$(id\o \D)(E)=(\D\o id)(E)=(1\o E)(E\o 1)=(E\o 1)(1\o E)$$
and $E\D(a)=\D(a)=\D(a)E$. In the third equality we use $m(E\rhd (r\o r'))=rr'$.

Finally, we have
\begin{align*}
E^B((a\o_D b)\cdot_D(a'\o_D b')) &= (E^B_1\rhd a)(E^B_2b_{(1)}\rhd a'\lhd S^{-1}(b_{(3)}))\o E^B_3b_{(2)}b'\\
                                  &= (E^B_1\rhd a\o_D E^B_2b )\cdot_D(a'\o_D b').
\end{align*}

So $R(A\o B)$ is a sub-algebra of $A\o B$. 
\end{proof}
{\bf Notation.} Now we will denote the algebra $R(A\o B)$ by $A\o_D B$.

\begin{proposition}\label{four}
For any $a,a', a''\in A, b, b', b''\in B$, the following  equalities hold:
	\begin{equation}\label{eq}
	\begin{split}
	&(a\o b)\cdot_D(a'\lhd \v_s(b')\o b'')=(a\o b)\cdot_D(a'\o \v_s(b')b'')\\
	&(a'\o \v_s(b')b'')\cdot_D(a\o b)=(a' \lhd \v_s(b')\o b'')\cdot_D(a\o b)
	\end{split}
	\end{equation}
and
	\begin{equation}\label{eq2}
	\begin{split}
	&(a\o \v'_t(a')\rhd b)\cdot_D (a''\o b')=(a\v'_t(a')\o  b)\cdot_D (a''\o b')\\
	&	(a''\o b')\cdot_D (a\o \v'_t(a')\rhd b)=(a''\o b')\cdot_D (a\v'_t(a')\o  b)
	\end{split}
	\end{equation}	
and 	
\begin{equation}\label{eq3}
\begin{split}
&(a\o b)\cdot_D( \v_s(b')\rhd a' \o b'')=(a\o b)\cdot_D(a'\o S(\v_s(b'))b'')\\
&( \v_s(b')\rhd a' \o b'')\cdot_D(a\o b) =(a'\o S(\v_s(b'))b'')\cdot_D(a\o b)
\end{split}
\end{equation}
and
\begin{equation}\label{eq4}
\begin{split}
&(a\o b\lhd \v_t(a'))\cdot_D( a'' \o b')=(a\v_t(a')\o b)\cdot_D(a''\o b')\\
&( a'' \o b')\cdot_D (a\o b\lhd \v_t(a'))=( a'' \o b')\cdot_D (a\v_t(a')\o b)
\end{split}
\end{equation}
\end{proposition}
\begin{proof}
 For any $a,a'\in A, b, b', b''\in B$, we have
\begin{align*}
(a\o b)\cdot_D(a'\lhd \v_s(b')\o b'') &= a(b_{(1)}\rhd a'\lhd \v_s(b')S^{-1}(b_{(3)}))\o b_{(2)}b''\\
                               &= a(b_{(1)}\rhd a'\lhd S^{-1}(b_{(3)}))\o b_{(2)}\v_s(b')b'' \\
                               &= (a\o b)\cdot_D(a'\o \v_s(b')b'').
\end{align*}
In the second equality we use $E(y\o 1)=E(1\o S(y))$ where $y$ is in the source algebra.
\begin{align*}
(a'\o \v_s(b')b'')\cdot_D(a\o b) &= a'((\v_s(b')b'')_{(1)}\rhd a \lhd S^{-1}((\v_s(b')b'')_{(3)}))\o (\v_s(b')b'')_{(2)}b\\
                               &= a'((E^B_1b''_{(1)})\rhd a\lhd S^{-1}(E^B_3b''_{(3)})S^{-1}(\v_s(b')))\o E^B_2b''_{(2)}b\\
                               &= a'((b''_{(1)})\rhd a\lhd S^{-1}(b''_{(3)})S^{-1}(\v_s(b')))\o b''_{(2)}b\\
                               &= (a'\lhd \v_s(b'))((b''_{(1)})\rhd a\lhd S^{-1}(b''_{(3)}))\o b''_{(2)}b\\
                               &= (a' \lhd \v_s(b')\o b'')\cdot_D(a\o b)
\end{align*}
In the fourth equality we use $a(a'\lhd x)=(a\lhd S(x))a'$ where $x$ belongs to the target algebra. Note that in the regular case,  $\v_s$ and $\v'_s$ have the same range. The proof of other  equations is similar.
\end{proof}

\begin{remark}\label{cover}
	
(1) If the product  to be non-degenerate, then we can cancel  $a\o b$ or $a''\o b'$ in above equalities. The non-degeneracy of the product   will be proved in the next proposition.\\
(2) We will introduce some notations.  For any $a\in A, b\in B$, we know that
$$R(b\rhd a\o b')=(b'_{(1)}b)\rhd a\o b'_{(2)}.$$
Obviously $b'_{(1)}$ is well covered.  For simplicity,  we write
\begin{eqnarray}
R( a\o b')=b'_{(1)}\rhd a\o b'_{(2)}.
\end{eqnarray}
 Similarly,  we can  consider the map $RR^{op,cop}_1$,
\begin{align*}
RR^{op,cop}_1(a\o a'\rhd b\lhd a'') &= R(a_{(2)}\o b_{(1)}\left\langle S^{-1}(S(a')a_{(1)}), b_{(2)}\right\rangle )\\
                                    &= R(a_{(2)}\o S^{-1}(S(a')a_{(1)})\rhd b\lhd a'')\\
                                    &= a_{(2)}\o S^{-1}(S(a')a_{(1)})\rhd b\lhd (a''a_{(3)})
\end{align*}
Here $a_{(1)}$ is covered by $a'$ and $a_{(3)}$ is covered by $a''$. Now we can denote it by
\begin{eqnarray}
RR^{op,cop}_1(a\o b)= a_{(2)}\o S^{-1}(a_{(1)})\rhd b \lhd a_{(3)}.
\end{eqnarray}
\end{remark}

\begin{proposition}
	The product $\cdot_D$ in $A\o_D B$ is non-degenerate.
\end{proposition}
\begin{proof}
	For any $\sum a_{(i)}\o b_{(i)}\in A\o_D B$, assume that
	$$(\sum_i a_{(i)}\o_D b_{(i)}) (p\o_D q)=0$$
	for all $p\in A, q\in B$. Now we have
	$$\sum_i a_{(i)}(b_{(1)}\rhd p\lhd S^{-1}_B(b_{(3)}))\o b_{(2)}q=0.$$
	Apply $\D$ and $S_B$ , then we get
	$$\sum_i a_{(i)}(b_{(1)}\rhd p\lhd S^{-1}_B(b_{(3)}))\o S_B(b_{(2)})q\o b_{(3)}q'=0 $$
	for all $p\in A, q, q'\in B$. Now write $p$ as $b'\rhd p'$ and replace $b'$ by $S_B(b_{(2)})q$.  We have
	$$\sum_i (E^B_1\rhd a_{(i)})(p\lhd S^{-1}_B(b_{(2)}))\o E^B_2b_{(1)}q=0$$
	for all $p\in A$ and $q\in B$. Now repeat the above procedure, we can get
	$$\sum_i E^B_{(1)}\rhd a_{(i)}\lhd S^{-1}(E^B_{(3)})p\o E^B_{(2)}b_{(i)}q=0$$
	for all  $p\in A$ and $q\in B$. We can cancel $p, q$, then by Proposition \ref{four} we have
	$$\sum_i E^B_1 \rhd a_{(i)}\o E^B_2b_{(i)}=0.$$
	On the other hand, suppose that
	$(p\o_D q)(\sum_i a_{(i)}\o_D b_{(i)})=0$ for all $p\in A, q\in B$.  As above, we can obtain
	$$\sum_i p(E^B_1\rhd a_{(i)}\lhd S_B^{-1}(E^B_{(3)}))\o qE^B_2b_{(i)}=0$$
	for all $p, q$. Finally we have $\sum_i E^B_1\rhd a_{(i)}\o E^B_2b_{(i)}=0$.
\end{proof}

Now by the non-degeneracy of the product,  we obtain 
\begin{eqnarray}\label{4}
\begin{cases}
	a\o x\rhd b= ax\o b\\
	a\o b\lhd x=aS(x)\o b
\end{cases}
\begin{cases}
	a\lhd y\o  b= a\o yb\\
	y\rhd a\o b=a\o  S(y)b
\end{cases}
\end{eqnarray}
in algebra $A\o_D B$, for any $a\in A, b\in B, x\in A_t, y\in B_s$.

In (multiplier) Hopf algebras case, we know that  the twist map $T$ is  bijective from $A\o B$ to $B\o A$(see Proposition 1.1 in \cite{De1} or Proposition 3.5 in \cite{DrVD}). Now in weak multiplier Hopf algebras case,  the map $T$ is not bijective any more.  Unsurprisingly, we will show that the map $T$ is generalized inversed.

\begin{proposition}
There exists $T^{-1}: A\o B\rightarrow B\o A$ such that 
$$TT^{-1}T=T, \quad \quad T^{-1}TT^{-1}=T^{-1}.$$	
And we call $T^{-1}$ the generalized inverse of $T$. Denote $\tau (A\o_D B)$ by $B\o_D A$, here $\tau$ is the flip map, then $T$ is a bijective map from $A\o_D B$ to $B\o_D A$.
\end{proposition}

\begin{proof}
 Define
$$T^{-1}=\tau R^{op,cop} R': A\o B\ra B\o A.$$
 We will use the notations  in  Remark \ref{cover}.  For any $a\in A, b\in B$, we have
 \begin{align*}
 TT^{-1}T(b\o a) &= TT^{-1}(a_{(2)}\o S^{-1}(a_{(1)})\rhd b\lhd a_{(3)})\\
 &=T(\v'_s(a_{(1)})\rhd b\lhd \v'_s(a_{(3)})\o a_{(2)})\\
 &= a_{(2)}\o S^{-1}(a_{(1)})\rhd b\lhd a_{(3)}\\
 &=T(b\o a).
 \end{align*}
 With a similar computation, we can get $T^{-1}TT^{-1}(a\o b)=T^{-1}(a\o b)$.
If we are working with the algebra $A\o_D B$, then we have 
\begin{align*}
T^{-1}T(b\o_D a) &= \tau(\v'_s(b_{(1)})\rhd a\lhd \v'_s(b_{(3)})\o_D b_{(2)})\\
                 &\overset{(\ref{eq})}{=} \tau(\v'_s(b_{(1)})\rhd a \o_D b_{(2)})\\
                 & = \tau(E^B_1\rhd a \o_D E^B_2b)\\
                 &= b\o_D a
\end{align*}
and
\begin{align*}
T T^{-1}(a\o_D b)  &= \v'_s(b_{(1)})\rhd a\lhd \v'_s(b_{(3)})\o_D b_{(2)}\\
                 &= a\o_D b.
\end{align*}
 Note that   we can replace $a\o_D b$ by $b'\rhd a\lhd b''\o a'\rhd b\lhd a''$ in the above calculation, then we do not need to worry about the coverings.
\end{proof}

Recall that since  the map $R$ can be extended  to multiplier algebra $M(A)\o M(B)$, then we can extend $T$ to the map from  $M(B)\o M(A)$ to $M(A)\o M(B)$.

\begin{proposition}
	Let $(A, B, \left\langle ,\right\rangle  )$ be a WMHA pairing, for any $a\in A, b\in B, x\in \v_s(A), y\in \v'_s(B)$, we have
	$T(b\o_D 1)=1\o_D b$ and $T(1\o_D a)=a\o_D 1$.
\end{proposition}
\begin{proof}
	\begin{equation*}
	T(1\o_D a)=E_1\rhd a\lhd S^{-1}(E_3)\o_D E_2= a\o_D 1.
	\end{equation*}
	If we write $a$ as $b\rhd a'\lhd b'$, then $E_1$ and $E_3$ can be covered. In the second equality we use formula (\ref{4}).  The  proof of $T(b\o_D 1)=1\o b$  is similar.
\end{proof}

Now we come to the main concept in this section.

\begin{definition}
The algebra $A\o_D B$ defined above is called the {\itshape Drinfeld double} of weak multiplier Hopf algebra pairing$(A, B, \left\langle ,\right\rangle  )$.
\end{definition}

 When $A$ and $B$ are multiplier Hopf algebras, we have $E(A\o B)=A\o B$, so the double defined above is coincided with the Definition 3.4 in \cite{DrVD}. When $A$ and $B$ are weak Hopf algebras, it is the case in \cite{B3} or \cite{NTV}.
 If $B=\widehat{A}$, the  Drinfeld double $A\o_D B$ will be denoted by $D(A)$.
\newline

 \noindent
 {\bf Notation.} Given a pairing $(A, B, \left\langle ,\right\rangle  )$, for any $a\in A, b\in B,$  we will use the  notation
 $$T(b\o a)= \sum a^i\o b^i\quad  \quad T^{-1}(a\o b)= \sum b_j\o a_j,$$
 and
 $$T(b^i\o a^i)=\sum a^{iI}\o b^{iI}\quad  \quad T^{-1}(a_j\o b_j)=\sum b_{jJ}\o a_{jJ}.$$
 We also omit the summation notation in the following.

 \begin{lemma}\label{T}
 	Given WMHA pairing $(A, B, \left\langle ,\right\rangle  )$, for any $a, a'\in A, b, b'\in B$, we have
 	
 	$(1)\quad T(bb'\o a)=T(b\o a^i)(1\o b'^i),$
 	
 	$(2)\quad T(b\o aa')=(a^i\o 1)T(b^i\o a').$
 \end{lemma}
 \begin{proof}
 	(1) \begin{align*}
 	T(bb'\o a) &=(bb')_{(1)}\rhd a\lhd S^{-1}({(bb')_{(3)}})\o (bb')_{(2)}\\
 	&= b_{(1)}\rhd a^i\lhd S^{-1}(b_{(3)})\o b_{(2)}b^{'i}\\
 	&= T(b\o a^i)(1\o b^{'i})
 	\end{align*}
 	(2) The proof is similar.
 \end{proof}
 As a consequence we have
 \begin{equation}\label{ijequ}
 T(b_j\o a_ja')=(a\o 1)T(b\o a') \quad \text{and} \quad T(b'b_j\o a_j)=T(b'\o a)(1\o b).
 \end{equation}

In the end of this subsection, we consider the $*$-case.

\begin{proposition}\label{*algebra}
Let $(A, B, \left\langle ,\right\rangle )$ be a weak multiplier Hopf  $*$-algebra pairing.  Then the Drinfeld double $D=A\o_D B$ is a $*$-algebra with the $*_D$  operation $$(a\o_D b)^{*_D}=T(b^*\o a^*).$$
\end{proposition}
\begin{proof}
	Recall that $T=R\circ R^{op,cop,'}\circ \tau$, then we have
	\begin{align*}
	((a\o_D b)^{*_D})^{*_D}&=RR^{op,cop,'}R'R^{op,cop}(a\o_D b)\\
	&= RR'R^{op,cop,'}R^{op,cop}(a\o_D b)\\
	&=E^B\rhd(a_{(2)}\o_D \v'_t(a_{(1)})\rhd b)\\
	&=a\o_D b.
	\end{align*}
	Now let us check the antimultiplicativity of $*_D$. For any $a,a'\in A, b, b'\in B$, we will use notation
	\begin{equation}\label{T-notation}
T(b\o a)=a^i\o b^i=a^{i'}\o b^{i'}=a^{I}\o b^{I}=a^{I'}\o b^{I'}=a^P\o b^P
	\end{equation}
	Then
	\begin{align*}
	((a\o_D b)(a'\o_D b'))^{*_D}&= (aa'^i\o b^ib')^{*_D}\\
	&= (a'^{i*}a^*)^I\o(b'^*b^{i*})^I\\
	&=a'^{i*II'}a^{*i'P}\o b'^{*I'P}b^{i*Ii'}\\
	&=a'^{*I'}a^{*i'P}\o b'^{*I'P}b^{*i'}.
	\end{align*}
	In the third equality we use Lemma \ref{T}, in the last equation we  use
	\begin{equation*}\label{*_DT}
	(T(b\o a))^{*_D}=a^*\o b^*.
	\end{equation*}
	On the other side, we compute
	\begin{align*}
	(a'\o_D b')^{*_D}(a\o_D b)^{*_D}&= (a'^{*i}\o_D b^{'*i})(a^{*i}\o b^{*i})\\
	&= a'^{*i}a^{*Ii'}\o b'^{*ii'}b^{*I}.
	\end{align*}
	Then we obtain $((a\o_D b)(a'\o_D b'))^{*_D}=(a'\o_D b')^{*_D}(a\o_D b)^{*_D}$.
\end{proof}

 \begin{lemma}
	Assume that $(A, B, \left\langle ,\right\rangle )$ is a WMHA pairing and $D$ is the Drinfeld double.  Let $x\in M(A), y\in M(B), X\in M(A\o A)$ and $Y\in M(B\o B)$. Then $F(m\o n)$ defined by
	$$
	F(x\o y)_l:= (x_l\o id)T(y_l\o id)T^{-1}
	$$
	$$
	F(x\o y)_r:= (id\o y_r)T(id\o x_r)T^{-1}
	$$
	is  a multiplier in $M(D)$, $F'(X\o Y)$ defined by
	$$
	F'(X\o Y)_l:= (X_l)_{13}(T\o T)(Y_l)_{13}(T^{-1}\o T^{-1})
	$$
	$$
	F'(X\o Y)_r:= (Y_r)_{24}(T\o T)(X_r)_{24}(T^{-1}\o T^{-1})
	$$
	is a multiplier in $M(D\o D)$.
\end{lemma}
\begin{proof}
	For any $a, p\in A, b, q\in B$, we still use (\ref{T-notation})
	and $$T^{-1}(a\o b)=b_j\o a_j=b_J\o a_J.$$
	Now we compute
	\begin{align*}
	(a\o_D b)(F(x\o y)_l(p\o_D q))&=(a\o_D b)(((x_l\o id)T(y_l\o id)T^{-1})(p\o_D q))\\
	&=(a\o_D b)(x_l(p_j^{\,\,\, i})\o_D y_l(q_j)^i)\\
	&=ax_l(p_j^{\,\,\, i})^{I}\o b^Iy_l(q_j)^i\\
	&=(a_Jx_l(p_j^{\,\,\, i}))^I\o (b_J)^Iy_l(q_j)^i\\
	&= (x_r(a_J)p_j^{\,\,\, i})^I\o (b_J)^Iy_l(q_j)^i\\
	&= x_r(a_J)^{I'}p_j^{\,\,\, I}\o (y_r(b_J^{\,\,\,I'})q_j)^I\\
	&= x_r(a_J)^Ip^i\o y_r(b_J^{\,\,\, I})^iq
	\end{align*}

In the fourth equation we use Lemma \ref{T}, in the fifth and sixth equations we use the property of multiplier. On the other side, we have
\begin{align*}
F(x\o y)_r(a\o_D b)(p\o_D q)&= (x_r(a_j)^i\o_D y_r(b_j^{\,\,\, i}))(p\o_D q)\\
&= x_r(a_j)^ip^I\o_D y_r(b_j^{\,\,\, i})^Iq
\end{align*}
Then we have proved that $F(x\o y)$ is a multiplier in $M(D)$. The proof of $F'(x\o y)$ is similar.
\end{proof}
Recall that if $A$ is a $*$-algebra,  then $M(A)$ is a $*$-algebra with $*$-operation $\rho^*=(\rho_r^*, \rho_l ^*)$ where $f^*(a)=f(a^*)^*$ for any $a\in A, f\in End_k(A)$, here $\rho=(\rho_l, \rho_r)$ is a multiplier in $M(A)$.

\begin{proposition}
	Let $(A, B, \left\langle ,\right\rangle )$ be a weak multiplier Hopf $*$-algebra pairing and $D$ be the Drinfeld double. Then the following four maps
	$$f_1: M(A)\rightarrow M(D), x\mapsto F(x\o 1)$$
	$$f_2: M(B)\rightarrow M(D), y\mapsto F(1\o y)$$
	$$f_3: M(A\o A)\rightarrow M(D\o D), X\mapsto G(X\o 1)$$
	$$f_4: M(B\o B)\rightarrow M(D\o D), Y\mapsto G(1\o Y)$$
	are $*$-algebra morphisms.
\end{proposition}
\begin{proof}
	First let us consider the proof of $f_1$.   Obviously, $f_1$ is a algebra embedding. Now we will show that $f_1$ is a $*$-algebra morphism.
	\begin{align*}
	(f_1(x)^*)_l(a\o_D b)&=(f_1(x)_r((a\o_D b)^*))^*\\
	&= *_DT(b^*\o x_r(a^*))\\
	&=(x_r(a^*))^*\o b=(x_l)^*(a)\o b\\
	&=f_1(x^*)_l(a\o_D b).
	\end{align*}
	In the third equation we use formula (\ref{*_DT}). The other cases can be proved in a similar way.
\end{proof}

\subsection{The weak multiplier Hopf algebra structure on $D$}

 Now we will show that there is a weak multiplier Hopf algebra structure on $D$.  

\begin{theorem}\label{main}
If $(A, B, \left\langle ,\right\rangle )$ is a WMHA  pairing, then $D$ is a regular weak multiplier Hopf algebra with the following data:
\begin{itemize}
\item  the coproduct $\D_D$ is defined as
$$
\D_D= F'(\D_A\o \D_{B^{cop}}): D\ra M(D\o D).
$$

\item the canonical idempotent $E_D$ is defined as
$$E_D(a\o_D b\o a'\o_D b')= E^A_1a\o_D E^B_2b\o E^A_2a'\o_D E^B_1b'.$$

\item the counit $\v_D$ is defined as $\v_D= D\ra \mathbb{C}: a\o_D b\mapsto \left\langle a, \v'_s(b)\right\rangle  $.

\item the antipode $S_D$ is defined as $S_D=T\circ\tau \circ(S_A\o S^{-1}_B)$.

\end{itemize}
If $A$ and $B$ are weak multiplier Hopf $*$-algebras, then the Drinfeld double is a weak multiplier Hopf $*$-algebra with the $*$-operation defined in Proposition \ref{*algebra}.

\end{theorem}
We will prove this theorem by several propositions.

\begin{proposition}
	$\D_D$ is coassociative and a homomorphism from $D$ to $M(D\o D)$.
\end{proposition}
\begin{proof}
First we will show that $\D_D$ is a comultiplication on $D$. According to the definition of $\D_D$,  the canonical maps $T^D_1, T^D_2$ on $D\o D$ are defined as
$$
T^D_1=(T^A_1)_{13}\circ(id\o id\o T)\circ(id\o T^{B^{cop}}_1\o id)\circ(id\o id \o T^{-1})
$$
$$
T^D_2=(T^{B^{cop}}_2)_{24}\circ(T\o id\o id)\circ(id\o T^{A}_2\o id)\circ(T^{-1}\o id \o id).
$$
Obviously, we have
$$T^D_1(d\o d')\in D\o D \quad \text{and} \quad T^D_2(d\o d')\in D\o D$$
 for any $d, d'\in D$.  From the  coassociativity of $\D_A$ and $\D_B$ we can obtain
 $$
 (id\o T^D_1)(T^D_2\o id)=(T^D_2\o id)(id\o T^D_1).
 $$
It means that $\D_D$ is coassociative. So $\D_D$ is a comultiplication.

Next we will show that  $\D_D$ is a homomorphism from $D$ to $M(D\o D)$.
\begin{align*}
& \D_D(a\o_D b)\D(a'\o_D b') \\
&= (a_{(1)}\o_D b_{(2)}\o a_{(2)}\o_D b_{(1)})(a'_{(1)}\o_D b'_{(2)}\o a'_{(2)}\o_D b'_{(1)})\\
&=a_{(1)}a'_{(2)}\o_D b_{(4)}b'_{(2)}\o a_{(2)}a'_{(4)}\o_D b_{(2)}b'_{(1)}\left\langle a'_{(1)}, S_B^{-1}(b_{(5)})\right\rangle \left\langle a'_{(3)}, \v'_s(b_{(3)}) \right\rangle  \left\langle a'_{(5)}, b_{(1)} \right\rangle \\
&= a_{(1)}a'_{(2)}\o_D b_{(4)}b'_{(2)}\o a_{(2)}a'_{(3)}\lhd \v'_s(b_{(3)})\o_D b_{(2)}b'_{(1)} \left\langle a'_{(1)}, S_B^{-1}(b_{(4)})\right\rangle  \left\langle a'_{(4)}, b_{(1)} \right\rangle \\
&\overset{(\ref{eq})}{=} a_{(1)}a'_{(2)}\o_D b_{(3)}b'_{(2)}\o a_{(2)}a'_{(3)}\o_D b_{(2)}b'_{(1)} \left\langle a'_{(1)}, S_B^{-1}(b_{(4)})\right\rangle  \left\langle a'_{(4)}, b_{(1)} \right\rangle.
\end{align*}
On the other hand we have
\begin{align*}
& \D_D((a\o_D b)(a'\o_D b'))\\
&= \D_D(a(b_{(1)}\rhd a'\lhd S^{-1}(b_{(3)}))\o_D b_{(2)}b')\\
&= a_{(1)}(a'_{(1)}\lhd S_B^{-1}(b_{(4)}))\o_D b_{(3)}b'_{(2)}\o a_{(2)}(b_{(1)}\rhd a'_{(2)})\o_D b_{(2)}b'_{(1)}\\
&= a_{(1)}a'_{(2)}\o_D b_{(3)}b'_{(2)}\o a_{(2)}a'_{(3)}\o_D b_{(2)}b'_{(1)} \left\langle a'_{(1)}, S_B^{-1}(b_{(4)})\right\rangle  \left\langle a'_{(4)}, b_{(1)} \right\rangle.
\end{align*}
In the second equality we use $\D_A(b\rhd a \lhd b')=a_{(1)}\lhd b' \o b\rhd a_{(2)}.$ Then we finish the proof. And the coverings will be explained in  following remark.
\end{proof}

\begin{remark}
Let us consider the range of map $T^D_1$. For any $a,a'\in A, b,n'\in B$, we have
\begin{align*}
T^D_1(a\o_D b\o a'\o_D b')&= a_{(1)}\o b_{(2)}\o (a_{(2)}\o 1)T(b_{(1)}b'_j\o a'_j)\\
&\overset{(4.7)}{=} a_{(1)}\o b_{(2)}\o (a_{(2)}\o 1)T(b_{(1)}\o a')(1\o b')\\
&=   a_{(1)}\o_D b_{(2)}\o ((a_{(2)}\o_D b_{(1)})(a'\o b')).
\end{align*}
The right side of first equality can be written as
$\sum a_{(1)}\o b_{(2)}\o a_{(2)}(a'_j)^i\o (b_{(1)}b'_j)^i,$
now $a_{(2)}$ is covered by $(a'_j)^i$ and $b_{(1)}$ is covered by $b'_j$. Thus, we can write $\D_D$ as
$$\D_D(a\o_D b)=\sum  a_{(1)}\o_D b_{(2)}\o a_{(2)}\o_D b_{(1)}.$$
This is a familiar formula is (weak) Hopf algebras theory. Moreover, since $\D_A$ and $\D_B$ are full, it is easy to get that $\D_D$ is full.

In  Proposition 4.14, if we multiply  with an element $d\in D$ (left or right) and write $a$ as $b\rhd a'\lhd b'$, then everything will be  covered. But  the proof will be too long. 
\end{remark}

\begin{lemma}\label{DeltaT}
	We have $\D_DT=(T\o T)(\tau \o \tau)\D_D\tau.$
\end{lemma}
\begin{proof}
	For any $a\in A, b\in B$,  we have
	\begin{align*}
	&(T\o T)(\tau \o \tau)\D_D\tau(b\o a) \\
	&=T(b_{(2)}\o a_{(1)})\o T(b_{(1)}\o a_{(2)})\\
	&= a_{(2)}\o_D b_{(4)}\o a_{(3)}\lhd \v'_t(b_{(3)})\o_D b_{(2)}\left\langle a_{(1)}, S^{-1}(b_{(5)}) \right\rangle
	\left\langle a_{(4)}, b_{(1)} \right\rangle\\
	&=	a_{(2)}\o_D b_{(3)}\o a_{(3)}\o_D b_{(2)}\left\langle a_{(1)}, S^{-1}(b_{(4)}) \right\rangle\left\langle a_{(4)}, b_{(1)} \right\rangle\\
	&=\D_DT(b\o a)
	\end{align*}
\end{proof}

\begin{proposition}\label{Delta*}
	Let $\left\langle A, B \right\rangle$ be a weak multiplier Hopf algebra $*$-pairing. Then the coproduct  is a $*$-algebra morphism.
\end{proposition}
\begin{proof}
	For any $a\o_D b\in D, c\o_D d\in D$,
	\begin{align*}
	&(*_D\o *_D)\D_D(a\o_D b)\\
	&= (T\o T)(*_A\o *_B\o *_A\o *_B)(\tau\o \tau)\D_D(a\o_D b)\\
	&= (T\o T)(\tau \o \tau)\D_D(*_A\o *_B)(a\o_D b)\\
	&=\D_DT\tau (*_A\o *_B)(a\o_D b)\\
	&=\D_D(*_D(a\o_D b))
	\end{align*}
	In the third equation we use Lemma \ref{DeltaT}.
\end{proof}

\begin{proposition}
	$\v_D$ is a counit.
\end{proposition}
\begin{proof}
For any $a,a'\in A, b, b'\in B$,  we have
\begin{align*}
&(id\o \v_D)T^D_2(a\o_D b\o a'\o_D b')\\
&= (id_D\o \v_D)(T(b_j\o a_ja'_{(1)})(1\o b'_{(2)})\o a'_{(2)}\o b'_{(2)})\\
&= (id_D\o \v_D)((b_{j(1)}\rhd a_i\lhd S^{-1}(b_{j((3))}))a'_{(2)}\o (S^{-1}(a'_{(1)})\rhd b_{i(2)}\lhd a'_{(3)})b'_{(2)}\o a'_{(4)}\o b'_{(1)})\\
&= (id_D\o \v_D)(aa'_{(2)}\o (S^{-1}(a'_{(1)})\rhd b_{i(2)}\lhd a'_{(3)})b'_{(2)}\o a'_{(4)}\o b'_{(1)})\\
&= (id_D\o \v_D)((a\o 1)T(b\o a'_{(1)})(1\o b'_{(2)})\o a'_{(2)}\o b'_{(1)})\\
&= (a\o 1)T(b\o a'_{(1)})(1\o b'_{(2)})\left\langle a'_{(2)}, \v'_s(b'_{(1)})\right\rangle  \\
&=  (a\o_D b)(\v'_s(b'_{(1)})\rhd a'\o_D b'_{(2)})\\
&\stackrel{(\ref{4})}{=} (a\o_D b)(a'\o_D b').
\end{align*}
Similarly, we can prove $(\v_D\o id)T^D_1(d\o d')=dd' $ for any $d, d'\in D$.
\end{proof}

Now let us consider the antipode $S_D$. First we have the following lemma.

\begin{lemma}
$S_D=(S_A\o S^{-1}_B)\circ \tau \circ T^{-1}.$
\end{lemma}
\begin{proof}
	By the definition of $S_D$, it is equivalent to show that
	$$(S_A\o S^{-1}_B)\circ \tau \circ T^{-1}= T\circ \tau \circ (S_A\o S_B^{-1}) .$$
	The formula  is equivalent with
	$$(S_A\o S_B^{-1})\circ T=\tau\circ  T^{-1}(S_A\o S_B^{-1})\circ  \tau .$$
	Recall the definition of $T$, then we need to show
	$$(S_A\o S_B^{-1})RR^{op,cop,'}=R^{op.cop}R'(S_A\o S_B^{-1}).$$
	From Proposition 3.29, we know that the equality above is true. Then we finish the proof.
\end{proof}

As a corollary we can get
$$S_D^2=S^2_A\o S^{-2}_B$$
and
$$S_D(a\o_D b)=S_A(a_j)\o_D S_B^{-1}(b_j)=S_A(a)^i\o_D S_B^{-1}(b)^i.$$
The two formulas will be used in the following.

\begin{proposition}
	The antipode $S_D$ defined in Theorem \ref{main} satisfies
	$$\sum d_{(1)}S_D(d_{(2)})d_{(3)}=d\quad  \text{and} \quad \sum S_D(d_{(1)})d_{(2)}S_D(d_{3})=S_D(d)$$
	for any $d\in D$.
\end{proposition}
\begin{proof}
For any $a\in A, b\in B,$ we compute
\begin{align*}
& \mu_D(\mu_D\o id)(id\o S_D\o id)(id\o \D_D)\D_D(a\o_D b) \\
 &= (a_{(1)}\o_D b_{(3)})S_D(a_{(2)}\o_D b_{(2)})(a_{(3)}\o_D b_{(1)})\\
 &= (a_{(1)}\o_D b_{(3)})(S_A(a_{(2)j})\o_D S_B^{-1}(b_{(2)j}))(a_{(3)}\o_D b_{(1)})\\
 &= (a_{(1)}S_A(a_{(2)j})^i\o_D (b_{(3)})^i S^{-1}_B(b_{(2)j}))(a_{(3)}\o_D b_{(1)})\\
 &= (a_{(1)}S_A(a_{(2)})^{Ii}  \o_D  (b_{(3)})^i S^{-1}_B(b_{(2)})^I) (a_{(3)} \o_D b_{(1)})\\
 &= a_{(1)}S_A(a_{(2)})^ia^I_{(3)}\o_D (\v'_s(b_{(2)}))^{iI}b_{(1)}\\
 &= a_{(1)}\v_s(a_{(2)})^i\o_D (\v'_s(b_{(2)}))^ib_{(1)}\\
 &\stackrel{(\ref{4})}{=} a\o_D b.
\end{align*}
Necessary coverings are explained in Lemma 4.8 and Remark 4.10. For the second formula we have
\begin{align*}
& S_D(a_{(1)}\o_D b_{(3)})(a_{(2)}\o_D b_{(2)})S_D(a_{(3)}\o b_{(1)}) \\
 &= (S_A(a_{(1)}^i)\o_D S_B^{-1}(b_{(3)})^i)(a_{(2)}\o_D b_{(2)})(S_A(a_{(3)J})\o_D S_B^{-1}(b_{(1)J}))\\
 &= (S_A(a_{(1)})^ia^I_{(2)}\o_D S_B^{-1}(b_{(3)})^{iI}b_{(2)})(S_A(a_{(3)J})\o_D S_B^{-1}(b_{(1)J}))\\
 &= (S_A(a_{(1)})a_{(2)})^iS_A(a_{3J})^{II'}\o_D S_B^{-1}(b_{(3)})^{iI'}b^I_{(2)}S_B^{-1}(b_{(1)J})\\
 &= (S_A(a_{(1)})a_{(2)})^iS_A(a_{3})^{i'I'}\o_D S_B^{-1}(b_{(3)})^{iI'}(b_{(2)}S_B^{-1}(b_{(1)}))^{i'}\\
 &\stackrel{(\ref{4})}{=} S_A(a)^i\o_D S_B^{-1}(b)^i\\
 &= S_D(a\o_D b).
\end{align*}
\end{proof}

\begin{proposition}
	For any $d,d'\in D$, we have
	$$T^D_1P^D_1(d\o d')=E_D(d\o d') \quad \text{and} \quad  T^D_2P^D_2(d\o d')=(d\o d')E_D.$$
	Moreover, the canonical idempotent $E_D$ satisfies
	$$(id\o D_D)(E_A)=(\D_D\o id)(E_D)=(E_D\o 1)(1\o E_D)=(1\o E_D)(E_D\o 1).$$
\end{proposition}

\begin{proof}
 We have
\begin{align*}
& T^D_1P^D_1(a\o_D b\o a'\o_D b')\\
&= T_1^D(a_{(1)}\o_D b_{(2)}\o S(a_{(2)}\o_D b_{(1)})(a'\o_D b'))\\
&= a_{(1)}\o_D b_{(3)} \o (a_{(2)}\o_D b_{(2)})(S(a_{(3)j})a'^i\o_D S^{-1}(b_{(1)j})^ib')\\
&= a_{(1)}\o_D b_{(3)} \o a_{(2)}(S(a_{(3)j})a'^i)^I\o_D b^I_{(2)}S^{-1}(b_{(1)j})^ib'\\
&= a_{(1)}\o_D b_{(3)} \o a_{(2)}(S(a_{(3)})a')^i\o_D (b_{(2)}S^{-1}(b_{(1)}))^ib'\\
&= a_{(1)}\o_D b_{(3)} \o a_{(2)}S(a_{(3)})a'\o_D b_{(2)}S^{-1}(b_{(1)})b'\\
&= E^A_1a\o_D E^B_2b\o E^A_2a'\o_D E^B_1b'\\
&= E_D(a\o_D b\o a'\o_D b')
\end{align*}
Similarly, we can get $T^D_2P^D_2(a\o_D b\o a'\o_D b')=(a\o_D b\o a'\o_D b')E_D$.

Finally, it is enough to prove
$$(\D_D\o id)(E)=(E_D\o 1)(1\o E_D).$$
This is easy to obtain because we have
$$(\D_A\o id)(E_A)=(E_A\o 1)(1\o E_A)\quad \text{and}\quad (\D_B\o id)(E_B)=(E_B\o 1)(1\o E_B).$$
\end{proof}

Now,  we have finished the proof of Theorem \ref{main}.  Next we will study the integral on  $D$. And we have the following result.

\begin{proposition}\label{int}
Assume that $(A, B, \left\langle ,\right\rangle )$ is a WMHA pairing. Let $\psi_A (\vp_A)$ and $\psi_B (\vp_B)$ be the right (left) integral on $A$ and $B$, respectively. Then $\psi_A\o \vp_B$ $ ( \varphi_A \o \psi_B)$ is a right (left) integral on $D$.
\end{proposition}
\begin{proof}
Let $a\in A, b\in B$, denote $(\psi_A\o \vp_B\o id)\D_D(a\o_D b)$ by $y$. Then it is enough to show
$\D_D(y)=(1\o y)E_D$.
\begin{align*}
& \D_D(y)\D_D(a'\o_D b')(1_D\o a''\o_D b'')\\
& =\D_D(y(a'\o_D b'))(1_D\o a''\o_D b'')\\
& = \D_D(\psi_A(a_{(1)})a_{(2)}a'^i\o_D \vp_B(b_{(2)})(b_{(1)})^ib')(1_D\o a''\o_D b'')\\
& =((\psi_A(a_{(1)})a_{(2)}(a'^i)_{(1)}\o_D (b_{(1)}\vp_B(b_{(2)}))^i_{(2)}b'_{(2)}\o a_{(3)}(a'^i)_{(2)}\o_D (b_{(1)}\vp_B(b_{(2)}))^i_{(1)}b'_{(1)}))\\
&\quad (1_D\o a''\o_D b'')\\
& =( (E^A_1(a^{'i})_{(1)}\o_D  ((b_{(1)}\vp_B(b_{(2)})^i))_{(2)}b'_{(2)} \o \psi_A(a_{(1)})a_{(2)}E^A_2(a'^i)_{(2)}\o  ((b_{(1)}\vp_B(b_{(2)})^i))_{(1)}b'_{(1)}\\
&\quad (1_D\o a''\o_D b'')\\
& = ((a'_{(1)}\o_D b'_{(2)}\o \psi_A(a_{(1)})a_{(2)}(a'_{(2)})^i\o_D (b_{(1)}\vp_B(b_{(2)}))^i b'_{(1)}))(1_D\o a''\o_D b'')\\
&=a'_{(1)}\o_D b'_{(2)}\o \psi_A(a_{(1)})a_{(2)}(a'_{(2)})^i(a'')^I\o_D ((b_{(1)}\vp_B(b_{(2)}))^i b'_{(1)}))^Ib''\\
&= a'_{(1)}\o_D b'_{(2)}\o \psi_A(a_{(1)})a_{(2)}(a'_{(2)})^i(a'')^{Ii'}\o_D (b_{(1)}\vp_B(b_{(2)}))^{ii'}(b'_{(1)})^Ib''
\end{align*}

On the other hand, we have
\begin{align*}
& (1\o y)E_D\D_D(a'\o_D b')(1\o a''\o_D b'')\\
&= (1\o y)(a'_{(1)}\o_D b'_{(2)}\o a_{(2)}a''^i\o_D (b'_{(1)})^ib'')\\
&= a'_{(1)}\o_D b'_{(2)}\o ((\psi_A\o \vp_B\o id)\D_D(a\o_D b)(1\o a_{(2)}a''^i\o_D (b'_{(1)})^ib''))\\
&= a'_{(1)}\o_D b'_{(2)}\o \psi_A(a_{(1)})a_{(2)}(a'_{(2)}a''^i)^I\o (\vp_B(b_{(2)})b_{(1)})^I(b'_{(1)})^ib''\\
&= a'_{(1)}\o_D b'_{(2)}\o \psi_A(a_{(1)})a_{(2)}(a'_{(2)})^I(a'')^{iI'}\o_D (\vp_B(b_{(2)})b_{(1)})^{II'}(b'_{(1)})^ib''
\end{align*}
So we have $\D_D(y)=(1\o y)E_D=E_D(1\o y)$, it means that $\psi_A\o \vp_B$ is a right integral on $D$. The proof of $\varphi_A\o \psi_B$ is similar.
\end{proof}

\begin{theorem}\label{Th1}
	Let $\left\langle A, B \right\rangle$ be a WMHA pairing. If $A$ and $B$ are algebraic quantum groupoids, then the related Drinfeld double $D$ is an algebraic quantum groupoid.
\end{theorem}
\begin{proof}
	 We have proved that there exists a WMHA structure on $D$. Now we only need to find a faithful set of integrals.
	 Assume the faithful set of right  integrals on $A$ is $\int_{r}^{A}$ and the faithful set of left integrals on $B$ is $\int_{l}^{B}$.    For any $a\o_D b \in D$, we have
	  \begin{align*}
	  \v^D_t(a\o_D b)&= (a_{(1)}\o_D b_{(2)})S_D(a_{(2)}\o_D b_{(1)})\\
	                            &= a_{(1)}S_A(a_{(2)}j)^i\o_D (b_{(2)})^iS^{-1}_B(b_{(1)j})\\
	                            &\stackrel{\ref{ijequ}}{=}  a_{(1)}S_A(a_{(2)})^i\o_D (b_{(2)}S^{-1}_B(b_{(1)}))^i\\
	                            &= \v_t (a)\o_D \v'_s (b).
	  \end{align*}
Similarly, we also have $\v^D_s(a\o_D b)=\v_s(a)\o_D \v'_t(b)$.
By  Proposition \ref{int} and \cite[Proposition 1.18]{VW4},  we know that the target algebra $\v_t(D)$ is spanned by  elements of the form
	  $$(f\o id)(E_D(d\o 1))\quad \text{or} \quad  (f\o id)((d\o 1)E_D)$$
	  where $f\in \int_{r}^{A}\o \int_{l}^{B}$ and $d\in D$. It means that $ \int_{r}^{A}\o \int_{l}^{B}$ is  a faithful set of integrals on $D$.  This completes the proof.
\end{proof}

Given a  weak multiplier Hopf ($*$-)algebra pairing $(A, B, \left\langle ,\right\rangle  )$, we have constructed the Drinfeld double $D$. We show that there exists a weak multiplier Hopf ($*$-)algebra structure on $D$. Moreover, if $A$ and $B$ are algebraic quantum groupoids, we also show that the  double is a algebraic quantum groupoid. Our result generalizes the classic Drinfeld double construction of Hopf algebras (\cite{Ma,M}). And it also covers the results in multiplier Hopf algebras (see \cite{DrVD}) and weak Hopf algebras (see \cite{B3,NTV}).

Next we study the basic example.

\begin{example}\label{EX1}
For any groupoid $G$, let $A$ be the algebra of the complex functions with finite support in $G$ and pointwise product. For any $p,q\in G$, the coproduct on $A$ is defined as

\begin{equation*}
\D(f)(p, q)=
\begin{cases}
f(pq), & \text{if $pq$ is defined;} \\
0, & \text{otherwise.}
\end{cases}
\end{equation*}
The counit is defined as
$$\v(f)=\sum f(e)$$
where the sum is taken over all the units of $G$. The canonical idempotent $E$ is defined as
$$
E(p, q)=
              \begin{cases}
                1, & \hbox{if $pq$ is defined;} \\
                0, & \hbox{otherwise.}
              \end{cases}
$$
The antipode is defined as $(S(f))(p)=f(p^{-1})$. Then $A$ is a regular weak multiplier Hopf algebra.

The dual of $A$ is the groupoid algebra $\mathbb{C}G$ and we denote it by $B$. Let $\lambda_p$ be the canonical embedding of $p\in G$ in $\mathbb{C}G$.   The product is defined as
$$
\l_p\l_q=
              \begin{cases}
                \l_{pq}, & \hbox{if $pq$ is defined;} \\
                0, & \hbox{otherwise.}
              \end{cases}
$$
 The coproduct is given as $\D_B(\l_p)=\l_p\o \l_p$. The counit $\v_B$ is defined as $\v(\l_p)=1$. $E_B$ equals to $\sum \l_e\o \l_e$ where $e$ is the unit in $G$. The antipode is defined as $S_B(\l_p)=\l_{p^{-1}}$.

Now we have a natural  pairing $\left\langle , \right\rangle  $ between $A$ and $B$ which is given as $\left\langle f,\l_p\right\rangle  =f(p)$, $f\in A, \l_p\in B$. Then we can  consider the Drinfeld double $D$.

First we consider the rang of map $R$. For any $f,g,h\in A, p, q, r\in G$, we have
$$R(f\o \l_p)=\sum f_{(1)}\left\langle f_{(2)}, \l_p\right\rangle  \o \l_p=\l_p\rhd f\o \l_p.$$
Since the action is unital, we have $R(A\o B)=A\o B=D$. Suppose that $pq$ is defined, then
\begin{align*}
& ((f\o \l_p)(g\o \l_q))(r\o h)\\
& =\left\langle S^{-1}(g_{(1)}), \l_p\right\rangle  \left\langle g_{(3)}, \l_p\right\rangle(f  g_{(2)}\o \l_p\l_q)(r\o h)\\
& =g_{(1)}(p^{-1})g_{(3)}(p)(fg_{(2)}\o \l_{pq})(r\o h)\\
&= f(r)g(p^{-1}rp)h(pq)
\end{align*}
It means that
\begin{align*}
 (f\o \l_p)(g\o \l_q)=\begin{cases}
 fg(p^{-1}\cdot p)\o \l_{pq}, &\hbox{$pq$ is defined}\\
 0, &\hbox{otherwise}
 \end{cases}
\end{align*}
 The coproduct on $D$ is defined as
\begin{align*}
& \D_D(f\o \l_p)(q\o g\o r\o h)\\
&= (f_{(1)}\o \l_p\o f_{(2)}\o \l_p)(q\o g\o r\o h)\\
&= f_{(1)}(q)g(p)f_{(2)}(r)h(p)\\
&=\begin{cases}
      f(qr)(gh)(p), & \hbox{if $qr$ is defined;} \\
     0, & \hbox{otherwise.}
    \end{cases}
\end{align*}
So the coproduct can be denoted by $\D(d)(x\o y)=d(xy)$, here $d\in D, x, y\in G\o A$ and the product $xy$ is the usual tensor product.

The counit is defined as
$$\v_D(f\o \l_p)=\left\langle f,\v_t(\l_p)\right\rangle  =f(t(p)).$$
$t(p)$ means the target of $p$ in the groupoid $G$.

Next let us consider the idempotent element $E_D$.
\begin{align*}
& E_D(f\o \l_p\o g\o \l_q)(r\o h\o r'\o h')\\
&=(E^A_1f)(r)(E^A_2g)(r')h(p)h'(q)\\
&=\begin{cases}
      f(r)g(r')h(p)h'(q), & \hbox{if $rr'$ is defined;} \\
      0, & \hbox{otherwise.}
    \end{cases}
\end{align*}
If $G$ is a group, then $E_D$ is $1\o 1$.

Finally, we have
\begin{align*}
& S_D(f\o \l_p)(r\o h)\\
&= T\circ \tau (S(f)\o S^{-1}(\l_p))(r\o h)\\
&= f_{(1)}(p)f_{(2)}(r)f_{(3)}(p^{-1})h(p^{-1})\\
&=     \begin{cases}
      f(prp^{-1})h(p^{-1}), & \hbox{if $pr$ and $rp^{-1}$ are defined;} \\
       0, & \hbox{otherwise.}
     \end{cases}
\end{align*}
It means that $S_D(f\o \l_p)=f(p\cdot p^{-1})\o \l_{p^{-1}}$.
\end{example}

\subsection{Yetter-Drinfeld modules}

First we need to recall the notation of complete module which is introduced by Van Daele in \cite{VD2}. Let $A$ be a regular weak multiplier Hopf algebras and $V$ a vector space. Let $(X, \lhd, \rhd)$ be a non-degenerate $A$-bimodule. Let $Z$ be the space of pair $(\lambda, \rho)$ of linear maps from $A$ to $X$ satisfying 
$$a\lhd \lambda(a')=\rho(a)\rhd a'$$
 for all $a, a'\in A$. For any $z=(\lambda, \rho)\in Z$, define
$$a\rhd z=(a\lambda(\cdot), \rho(\cdot a))\quad \quad z\lhd a=(\lambda(a\cdot), \rho(\cdot)a),$$
then $Z$ is an $A$-bimodule with the submodule $X$. Now we call $Z$ the complete module of $A$ and denote by $M_0(X)$.

Note that the theory of comodules for weak multiplier bialgebras was studied in \cite{B2}. Now we will take a slightly different approach to  comodule theory through the complete module.

\begin{definition}
	Let $A$ be a weak multiplier Hopf algebra. A right comodule  for $A$ is a vector space with a comodule map $\delta: V\rightarrow M_0(V\o A)$ such that 
		$$(id\o \D)\delta=(\delta\o id)\delta.$$
	The right comodule is called {\itshape full} if 
	$$\{(id\o x)\delta(v)| \forall x\in \hat{A}, v\in V\}=V. $$

\end{definition}

As in \cite{VD2},	we use notation $\delta(v)=\sum v_{(0)}\o v_{(1)}$.   And the left comodule can be defined in a similarly way. When $V$ is a left comodule, we use $\delta(v)=\sum v_{(-1)}\o v_{(0)}$. For detail coverings about the Sweedler notation, we refer to Proposition 1.10 or Example 2.10 in \cite{VD2}.

\begin{proposition}\label{co-module}
	There is a one to one correspondence between left  $\hat{A}$-module and right $A$-comodule, via
		$$x\rhd v= (id\o x)\delta(v),  \quad \quad x\in \hat{A}, v\in V.$$  
	And the module is unital if and only if the comodule is full.
\end{proposition}
\begin{proof}
The proof is easy.
\end{proof}

Comodules are usually required to be counital as well, and we know that  for unital weak bialgebra, B\"{o}hm  showed that $\delta$ is counital if and only if $\delta$ is left full in paper \cite{B2}. For general weak multiplier Hopf algebras, we have following result.

\begin{proposition}
	Let $A$ be an algebraic quantum groupoid with $\varepsilon\in \hat{A}$. $(V, \delta)$ is a right $A$-comodule, then $\delta$ is left full if and only if $\delta$ is counital.
\end{proposition}
\begin{proof}
	Assume that $\delta$ is left full. Since $\varepsilon\in \hat{A}$, by Proposition \ref{co-module}, the left $\hat{A}$-module is unital, then $\delta$ is left full. Conversely, we have $(id\o \varepsilon)\delta(v)=v$, so $\delta$ is left full.
\end{proof}

 We call the map 
$$\rho: V\o A\rightarrow V\o A, \quad v\o a \mapsto \delta(v)(1\o a)$$
the corepresentation of $V$ on $A$. Define
$$\rho':  V\o A\rightarrow V\o A, \quad v\o a\mapsto v_0\o S(v_1)a,$$
then we have $\rho\rho'\rho=\rho$ and $\rho'\rho\rho'=\rho'$. This means that $\rho'$ is the generalized inverse of $\rho$. Note that the corepresentation theory for weak multiplier  bialgebras was developed in \cite{B2}.

Since the left action of $A$ on $V$ is unital, we can define map 
$$\Pi: A\o V\rightarrow A\o V, \quad a\o v\mapsto E_1a\o E_2\rhd v.$$
 The map $\Pi$ also appeared in   the definition of Yetter-Drinfeld module over weak multiplier bialgebra, see Definition 3.3 in \cite{B3}.

\begin{definition}
	Let $A$ be a regular weak multiplier Hopf algebras and $V$ a vector space. Then we call $V$ a left-left Yetter-Drinfeld module over $A$ if the following hold:
	\begin{itemize}
		\item $(V, \rhd)$ is a unital left $A$-module;
		\item $(V, \delta)$ is a left full $A$-comodule;
		\item  For any $a,a'\in A, v\in V$, $V$ satisfies the compatibility condition
		\begin{equation}\label{YD-1}
		\sum (a_{(1)}\rhd v)_{(-1)}a_{(2)}a'\o (a_{(1)}\rhd v)_{(0)}=a_{(1)}v_{(-1)}a'\o a_{(2)}\rhd v_{(0)}
		\end{equation}
		\item   We require
		\begin{equation}\label{YD-2}
		\Pi \delta(v)=\delta(v).
		\end{equation}
	\end{itemize} 
\end{definition}

Using Sweedler notation, the formula (\ref{YD-2}) can be denoted by 
$$\sum E_1v_{(-1)}\o E_2\rhd v_{(0)}=v_{(-1)}\o v_{(0)}.$$
If we multiply with an extra element of $A$ in the first factor, we get $v_{(-1)}$ covered. More equivalent descriptions of  formula (\ref{YD-2}) can be found in \cite{B3}.

In the next proposition $D(A)$ means the algebra $\hat{A}\o_D A$ with product 
$$(\omega\o_D a)(\eta\o_D b)=\sum \eta_{(2)}\omega \o_D a_{(2)}b\left\langle S(a_{(1)}), \eta_{(1)} \right\rangle \left\langle a_{(3)}, \eta_{(3)}\right\rangle,$$
where $a, b\in A, \omega,\eta\in \hat{A}$.

\begin{proposition}
	Let $A$ be  a  algebraic quantum groupoid, then the left-left Yetter-Drinfeld $A$-module can be identified with the left module over the Drinfeld double $D(A)$.
\end{proposition}
\begin{proof}
   Assume that $V$ is a left-left Yetter-Drinfeld module. For any $a,b \in A, \omega, \eta\in \hat{A}$, by Proposition \ref{co-module},   the action of $D(A)$ on $V$ is given as
 \begin{equation}
(\omega\o_D a)\rhd v=\left\langle (a\rhd v)_{(-1)}, \omega \right\rangle (a\rhd v)_{(0)}= (a\rhd v)\lhd \omega.
\end{equation}
Now we check  that $V$ is a left $D(A)$-module. Indeed, 
 \begin{align*}
( \omega\o_D a)\rhd(( \eta\o_D b)\rhd v )&=(a\o_D \omega)\rhd ((b\rhd v)\lhd \eta)\\
&= (a\rhd ((b\rhd v))\lhd \eta)\lhd \omega
 \end{align*}
and
 \begin{align*}
 ((\omega\o_D a)(\eta\o_D b))\rhd v &=(\eta_{(2)}\omega \o_D a_{(2)}b)\left\langle S(a_{(1)}), \eta_{(1)} \right\rangle \left\langle a_{(3)}, \eta_{(3)}\right\rangle\rhd v\\
 &\stackrel{(\ref{YD-1})}{=} \left\langle S(a_{(1)})a_{(2)}(b\rhd v)_{(-1)}), \eta\right\rangle((a_{(3)}\rhd (b\rhd v)_{(0)})\lhd \omega)\\
 &\stackrel{(\ref{YD-2})}{=}  \left\langle (b\rhd v)_{(-1)}), \eta\right\rangle((a\rhd (b\rhd v)_{(0)})\lhd \omega)\\
 &=(a\rhd ((b\rhd v))\lhd \eta)\lhd \omega.
 \end{align*}
 Obviously, the module is unital.  
\end{proof}

In the end of this section, we consider the relation between  Drinfeld double and the smash product. The smash product of WMHA was studied in \cite{ZhouW}.
Let $A$ be an algebraic quantum groupoid. We now consider the WMHA pairing $\left\langle \widehat{A},A \right\rangle$. From previous sections we can construct the Drinfeld double algebra $D(A)=\widehat{A}\o_D A$. Let $D(A)_0=\widehat{A}\o A$ be the algebra with product
$$(x\o a)\cdot_0(y\o b)=(x\o_D y_{(1)}\rhd a)(y_{(2)}\o_D b),$$
$x, y\in \widehat{A}, a,b \in A$.
Note that the action $\rhd$ is unital, then $y_{(1)}$ can be covered. The product $(x\o_D y_{(1)}\rhd a)(y_{(2)}\o_D b)$ is taken from $D(A)$. The new multiplication defined above was first appeared in \cite{Lu} where $A$ is a Hopf algebra.
Usually we have to show that the product is associative and non-degenerate, now the following result will imply these properties.

\begin{proposition}
	With the notions above, we have
	$$D(A)_0\cong \widehat{A}\# A.$$
\end{proposition}
\begin{proof}
	For any $x, y\in \widehat{A}, a,b \in A$, we compute
	\begin{align*}
	(x\o a)\cdot_0(y\o b)
	&= (x\o_D a_{(1)})(y_{(2)}\o_D b)\left\langle y_{(1)}, a_{(2)}\right\rangle\\
	&= xy_{(3)}\o_Da_{(2)}b\left\langle y_{(1)}, a_{(4)}\right\rangle \left\langle y_{(2)}, S^{-1}(a_{(3)}) \right\rangle \left\langle y_{(4)}, a_{(1)}\right\rangle\\
	&= xy_{(2)}\left\langle y_{(3)}, a_{(1)}\right\rangle \left\langle y_{(1)}, \v'_s(a_{(3)})\right\rangle \o_D a_{(2)}b\\
	&= x(E_1a_{(1)}\rhd y\lhd S^{-1}(a_{(3)}))\o_D E_2a_{(2)}b\\
	&= x(a_{(1)}\rhd y)\o_D a_{(2)}b\\
	&= (x\#a)(y\#b).
	\end{align*}
\end{proof}

\section{Quasitriangular weak multiplier Hopf algebra}\label{s5}
\subsection{main definition}

\begin{definition}
	Let $A$ be a regular WMHA.  The element $R\in M(A\o A)$ is called {\itshape $E$-inverse} if there exists $\bar{R}\in M(A\o A)$ satisfies
	$$R\bar{R}=E^{cop}, \quad \bar{R}R=E$$
	and
	$$RE=R, \quad \bar{R}E^{cop}=\bar{R}.$$
 We call $\bar{R}$ the $E$-inverse element of $R$.
\end{definition}

 Immediately, we  have
$$R\bar{R}R=R  \quad \text{and} \quad   \bar{R}R\bar{R}=\bar{R}.$$
Now  $\bar{R}$ is the {\itshape generalized inverse} of $R$.  If $A$ is Hopf algebra, the $E$-inverse element is just the usual invertible element.

\begin{proposition}
 The element	$\bar{R}$   is  unique.
\end{proposition}
\begin{proof}
 Assume $R'$ is  another $E$-inverse element, we have
	$$\bar{R}=\bar{R}R\bar{R}=E\bar{R}=R'R\bar{R}=R'E^{cop}=R'RR'=R'.$$
\end{proof}

\begin{definition}\label{QT}
	A quasitriangular weak multiplier Hopf algebra is a pair $(A, R)$ where $A$ is a regular weak multiplier Hopf algebra and $R\in M(A\o A)$ is an $E$-inverse element  satisfying the following conditions: \\
	(1) $(\D\o id)(R)=R^{13}R^{23}$,\\
	(2) $(id\o \D)(R)=R^{13}R^{12}$,\\
	(3) $\D^{cop}(a)R=R\D(a)  \text{ for all }  a\in A$.
\end{definition}

\begin{remark}
	Note that $(\D\o id)$ and $(id\o \D)$ can be extended to $M(A\o A)$, then the formulas  (1) and (2) in the above definition are meaningful.  Here $R^{12}=R\o 1, R^{23}=1\o R$, etc. as usual.	Similarly,  the third condition holds for any $m\in M(A)$, i.e. we have $\D^{cop}(m)R=R\D(m)$.
\end{remark}

In this section, QT-WMHA stands for quasitriangular weak multiplier Hopf algebra.

\begin{proposition}
	Let $(A, R)$ be a QT-WMHA, then $R$ satisfies the Yang-Baxter equation
	$$R^{23}R^{13}R^{13}=R^{12}R^{13}R^{23}.$$
\end{proposition}
\begin{proof}
First  we have
	$$(id\o \D^{cop})(R)(1\o R)=(1\o R)(id\o \D)(R).$$
	Now by  condition (1) and (2), the proof is completed.
\end{proof}

In the following we suppose that  $R(a\o 1), (a\o 1)R, R(1\o b) $ and $(1\o b)R$ belong to $A\o A$ for all $a, b\in A$, then we denote $R(a\o 1)=R_1a\o R_2\in A\o A$, similarly for other expressions.

\begin{proposition}
	We have $(\v\o id)(R)=(id\o \v)(R)=1.$
\end{proposition}
\begin{proof}
	By the above assumption, we can define two multipliers $(\v\o id)(R)$ and $(id\o \v)(R)$  in $M(A)$. Denote $(\v\o id)(R)=X, (id\o \v)(R)=Y$. First apply $\v\o id \o id$ and $id\o \v\o id$ to condition (1) of Definition \ref{QT}, then we have
    $$(1\o X)R=R=R(1\o X). $$
	Similarly, from condition (2) of Definition \ref{QT}  we can get
	$$  (Y\o 1)R=R=R(Y\o 1).$$
Next, we can multiply above equations by $\bar{R}$ and then we obtain
\begin{eqnarray}
E(1\o X)=E=E(Y\o 1)\\
(1\o X)E^{cop}=E^{cop}=(Y\o 1)E^{cop}.
\end{eqnarray}
Recall the result in \cite[Lemma 1,1]{VW3}, we know  that if $a\in A$ and $E(a\o 1)=0$, then $a=0$. Now by equation (5.1), we have $X=Y=1$.
\end{proof}

\begin{remark}
	 When $A$ is Hopf algebra, we have $E=1\o 1$ and the proof is simple,
	  see  \cite[Lemma 2.1.2]{Ma}. And  $E=1\o 1$ is also true for multiplier Hopf algebras. For weak Hopf algebras, we still have the same result since  the coproduct is full. 
\end{remark}

\begin{proposition}
	Let $(A, E, R)$ be a QT-WMHA. For any $x\in A_t, y\in A_s$, we have
	\begin{eqnarray}
	(1\o x)R=R(x\o 1)     \label{r1}\\
	(x\o 1)R=(1\o S(x))R \label{r2}\\
    R(1\o x)=R(S(x)\o 1)	\label{r3}\\
	R(y\o 1)=R(1\o S(y))  \label{r4} \\
	(y\o 1)R=R(1\o y)    \label{r5}\\
	(1\o y)R=(S(y)\o 1)R \label{r6}
	\end{eqnarray}
\end{proposition}
\begin{proof}
Recall that we have Lemma \ref{ASAT}.

(5.3): $R(x\o 1)=RE(x\o 1)=R\D(x)=\D^{cop}(x)R=(1\o x)R$.

(5.4): $(x\o 1)R=(x\o 1)E^{cop}R=(1\o S(x))E^{cop}R=(1\o S(x))R$.

(5.5) follows from $(1\o x)E=(S(x)\o 1)E$ and $E^{cop}R=R$.

The proof of (\ref{r4}-\ref{r6}) is similar.
\end{proof}

\begin{proposition}
	Let $(A, E, R, S)$ be a QT-WMHA. We have
	\begin{equation}
	(\v_s\o id)(R)=E \quad \text{ and } \quad (id\o \v'_s)(R)=E^{cop}
	\end{equation}
	\begin{equation}
	(\v_t\o id)(R)=E^{cop} \quad \text{and} \quad (id\o \v'_t)(R)=E
	\end{equation}
	\begin{equation}
	(id\o \v_s)(R)=(S\o id)E^{cop} \quad\text{and}\quad  (id\o \v_t)(R)=(S\o id)E
	\end{equation}
	\begin{equation}\label{12}
	(S\o id)(R)=(id\o S^{-1})(R)=\bar{R} \quad \text{and} \quad (S\o S)(R)=R
	\end{equation}
\end{proposition}
\begin{proof}
	First we need to show that the left sides of these equalities are well defined. Denote $(\v_s\o id)(R) $ by $X$, define
	$$X(a\o b)=\v_s(R_1)a\o R_2b,   \quad (a\o b)X=a\v_s(R_1)\o bR_2.$$
	Since  $R(1\o b)$ and $(1\o b)R$ are in $A\o A$, the definition is meaningful. It is easy to check that $X\in M(A\o A)$. Similarly, we can define other  multipliers.  Now we give the proof.
$$	(\v_s\o id)(R) = E_1\v(R_1E_2)\o R_2\stackrel{(\ref{r1})}{=}E((\v\o id)R\o 1)=E$$
	$$	(\v_t\o id)(R)=\v(E_1R_1)E_2\o R_2\stackrel{(\ref{r5})}{=}\v(R_1)E_2\o R_2E_1=E^{cop}	$$
    $$(id\o \v'_s)(R)=R_1\o E_1\v(E_2R_2)=R_1E_2\o E_1\v(R_2)=E^{cop} $$
    $$(id\o \v'_t)(R)=R_1\o \v(R_2E_1)E_2=E_1R_1\o \v(R_2)E_2=E$$
    $$(id\o \v_s)(R)=R_1\o E_1\v(R_2E_2) \stackrel{(\ref{r3})}{=} S(E_2)R_1\o E_1\v(R_2)=(S\o id)E^{cop}$$
    $$(id\o \v_t)(R)=R_1\o \v(E_1R_2)E_2\stackrel{(\ref{r6})}{=}S(E_1)R_1\o \v(R_2)E_2=(S\o id)E$$

    For (\ref{12}), apply $\mu_{12}(S\o id\o id)$ and $\mu_{12}(id \o S\o id)$  to condition (1) of Definition \ref{QT} respectively, then we have
    $$((S\o id)(R))R=E \quad \text{and} \quad R((S\o id)(R))=E^{cop}.$$
    By  uniqueness  we have $(S\o id)(R)=\bar{R}$.
    For  $(id\o S^{-1})(R)=\bar{R}$, the proof is similar.

     The last equation is  equivalent with $(id\o S)(\bar{R})=R$ which  can be obtained by applying $id\o S$ to  $(id\o S^{-1})(R)=\bar{R}$.  
\end{proof}

	 Now let us consider the multiplier $(\D\o id)(\bar{R})$. Denote $(\D\o id)(\bar{R})$  by $X$, $(\D\o id)(R)$ by $M$, $(\D\o id)(E)$ by $N$ and $(\D\o id)(E^{cop})$ by $Q$. Then we have
	 \begin{equation}\label{XMNQ}
	 \begin{split}
	 XM=N, NX=X, MN=M\\
	 MX=Q, XQ=X, QM=M.
	 \end{split}
	 \end{equation}
      From these equations we have
      \begin{equation}\label{XM}
      XMX=X \quad \text{and} \quad MXM=M.
      \end{equation}
       Next we will show that  $X=\bar{R}^{23}\bar{R}^{13}$  is a solution of equations (\ref{XMNQ}). Recall that $M=R^{13}R^{23}, N=E^{12}E^{23}$.
      \begin{align*}
      XM=\bar{R}^{23}\bar{R}^{13}R^{13}R^{23}&=\bar{R}^{23}E^{13}R^{23}\\
      &\stackrel{(\ref{r1})}{=}\bar{R}^{23}(E_1\o R_1E_2\o R_2)\\
      &= E_1\o E_2e_1\o e_2=N
      \end{align*}
      \begin{align*}
      PX=XMX&=\bar{R}^{23}\bar{R}^{13}R^{13}R^{23}\bar{R}^{23}\bar{R}^{13}\\
      &= \bar{R}^{23}E^{13}(E^{cop})^{23}\bar{R}^{13}\\
      &=\bar{R}^{23}(E^{cop})^{23}E^{13}\bar{R}^{13}\\
      &= \bar{R}^{23}\bar{R}^{13}=X
      \end{align*}
      The verification of $MX=Q$ and $XQ=X$ is similar. Finally, assume that $X'$ is another solution, then we have
      $$X'=NX'=XMX'=XQ=X.$$
      This give us the following proposition.
\begin{proposition}
	$(\D\o id)(\bar{R})=\bar{R}^{23}\bar{R}^{13}$. Moreover, $\bar{R}^{23}\bar{R}^{13}$ is the unique generalized inverse of $R^{13}R^{23}$. Similarly, $(id\o \D)(\bar{R})=\bar{R}^{12}\bar{R}^{13}$, and $\bar{R}^{12}\bar{R}^{13}$ is the unique generalized inverse of $R^{12}R^{13}$.
\end{proposition}

	Given a QT-WMHA $(A, S, R)$. For any $a\in A$, define $u\in M(A)$ by
	$$ua=\sum S(R_2)R_1a$$
	$$au=uS^{-2}(a).$$
	Now we have the following result.
\begin{proposition}
$u$ is an invertible  element in $M(A)$  obeying
$$\D(u)=\bar{R}\bar{R}^{21}(u\o u)=(u\o u)\bar{R}\bar{R}^{21}, \quad  \quad S^2(a)=uau^{-1}$$
for any $a\in A$. 
\end{proposition}
\begin{proof}
	Obviously, $u$ is a left multiplier. For any $a, b, c \in A$, we have
	$$R_1a_{1}b\o R_2a_{(2)}\o ca_{(3)}=a_{(2)}R_1b\o a_{(1)}R_2\o ca_{(3)},$$
	then
	$$S^2(ca_{(3)})S(R_2a_{(2)})R_1a_{(1)}b=S^2(ca_{(3)})S(a_{(1)}R_2)a_{(2)}R_1b.$$
	By formula (\ref{r5}) and $RE=R$, we have
	$$S^2(c)uab=S^2(ac)ub.$$
	This means that $u$ is a multiplier of $A$.

	For any $a\in A$, define left multiplier $w$ of $A$ as follows:
	   $$wa=\sum S^{-1}(\bar{R}_2)\bar{R}_1a$$
	   $$aw=\sum wS^2(a).$$
	 Note that  we have
$\bar{R}(a\o 1), (a\o 1)\bar{R}, \bar{R}(1\o b) $ and $(1\o b)\bar{R}$ belong to $A\o A$ for all $a, b\in A$. So the above definition is meaningful.  In a similar way we can show that $w$ is a multiplier. Then we have
 \begin{align*}
 (uw)a
 &= u(S^{-1}(\bar{R}_2)\bar{R}_1a)=S^{-1}(\bar{R}_2)u\bar{R}_1a\\
 &= S^{-1}(E_1)E_2a=a,
 \end{align*}
so $uw=1\in M(A)$. 

In the following the letters $r, Q, q, T, t$ are different copies of $R$. We compute
\begin{align*}
\D(u)&= S(R_{2(2)})R_{1(1)}\o S(R_{2(1)})R_{1(2)}\\
&= S(R_{2(2)}r_{2(2)})R_1\o S(R_{2(1)}r_{2(1)})r_1\\
&= S(T_2Q_2)T_1t_1\o S(t_2q_2)Q_1q_1\\
&=S(\bar{R}_1Q_2)u\o S(\bar{R}_2q_2)Q_1q_1\\
&=\bar{R}_1S(Q_2)u\o \bar{R}_2uQ_1	  \\       
&= \bar{R}\bar{R}^{21}(u\o u),
\end{align*}
 In the third equality we use the conditions in Definition \ref{QT}. In the fourth equality we use formula (\ref{12}).  The fifth equality  uses 
$$Q_1q_1\o \bar{R}_1Q_2\o \bar{R}_2q_2=q_1Q_1\o Q_2\bar{R}_1\o q_2\bar{R}_2$$
which can be derived from the Yang-Baxter equation. To see the commute relation we deduce  from $au=uS^{-2}(a)$ and $(S\o S)(\bar{R})=R$. By the definition of $u$, we get the last equation easily.

\end{proof}

Denote $v=S(u)$,  we have
$$\D(v)=(v\o v)\bar{R}\bar{R}^{21}=\bar{R}\bar{R}^{21}(v\o v) \quad  \text{and}\quad S^{-2}(a)=vav^{-1},$$
for any $a\in A$. The proof is similar. 

\begin{proposition}\label{S^4}
	Let $(A, R)$ be a QT-WMHA with the element $u, v$ defined as above, then
	 
	$(1)$ $uv=vu$ and $uv$ is central;
	
	$(2)$ $h=uv^{-1}$ is a grouplike element in $M(A)$  and $$S^4(a)=hah^{-1}$$
	 for any $a\in A$.
\end{proposition}
\begin{proof}
	 (1) Since $S^2(a)=uau^{-1}$ holds for all $a\in A$, it can be extended to $M(A)$. Now we have $S^2(u)=u$ and $S^2(v)=v$. So $uv=S^2(v)u=vu$. Take any $a\in A$, 
	 then we get $uva=uS^{-2}(a)v=auv.$
	 
	 (2) Now apply $S$ to  $S^2(a)=uau^{-1}$,  then we have $$S^3(a)=S(u^{-1})S(a)S(u).$$
	 Since $S(A)=A,$ it means
	 $$S^2(a)=S(u^{-1})aS(u).$$

	 Note that we have $S(u^{-1})=S(u)^{-1}=(R_1S(R_2))^{-1}$. 
	 We compute
	 \begin{align*}
	 \D(S(u)) &= (\tau(S\o S))(\bar{R}\tau(\bar{R})(u\o u))\\
	             &= (S(u)\o S(u))\tau(\tau(\bar{R})\bar{R})\\
	             &= (S(u)\o S(u))\bar{R}\tau(\bar{R}).
	 \end{align*}
	 Then $$\D(S(u^{-1}))=(S(u^{-1})\o S(u^{-1}))R^{21}R.$$ 
	 And now we can obtain $\D(uS^{-1}(u))=uS^{-1}(u)\o uS^{-1}(u)$. 
	 
	 Finally,
	 $$S^4(a)=S^2(S(u^{-1})aS(u))=uS(u^{-1})aS(u)u^{-1}=hah^{-1}.$$
\end{proof}

Let $(A, R)$ be a QT-WMHA. For convenience, denote $W=R^{21}R$. Note that for any $f\in A'$, we can define the multiplier $(f\o id)(Q)\in M(A)$. 

\begin{proposition}
	With the notation as above,  Define map 
	$$F: \hat{A}\rightarrow M(A), \varphi\mapsto (\varphi \o id)(W).$$
	Then the image of $F$ lies in the centralizer of $A_s$.
\end{proposition}
\begin{proof}
	Take any $a\in A, y\in A_s$, we have 
	\begin{align*}
	F(\varphi)ya=[(\varphi\o id)(W)y]a&=(\varphi \o id)[W(1\o ya)]\\
	&= \varphi(r_2R_1)r_1R_2ya\\
	&\stackrel{(5.7)}{=}\varphi(r_2yR_1)r_1R_2a\\
	&=yF(\varphi)a.
	\end{align*}
	Therefore we get $F(\varphi)y=yF(\varphi)$.
\end{proof}

\begin{definition}\label{factorisable}
	Let $(A, R)$ be a QT-WMHA. $A$ is called {\it factorisable} if the map defined above is surjective. Equivalently, the  map
	$$F': \hat{A}\rightarrow M(A), \varphi\mapsto (id \o \varphi)(W)$$
	is surjective.
\end{definition}
From formula (\ref{12}), we can get 
 $$(S\o S)(W)=W^{12}.$$
 With the fact that the  antipode $S$ is invertible, the equivalence in Definition \ref{factorisable}  can be derived. 
\newline

Next we give the final  result of this subsection. Recall that the left adjoint action of $A$ on $A$ is 
$$h\cdot g=\sum h_{(1)}gS(h_{(2)}), \quad  \forall g, h\in A$$
 Note that  the adjoint action is not unital,  and we have to consider a sub-algebra $A_0$ of $A$, then $A_0$ is an $A$-module algebra, see  Proposition 3.11 in \cite{ZhouW}. The action will be extended to tensor square in usual way.

\begin{proposition}
	The element $(S\o id)(W)$ in $M(A\o A)$ is invariant under the adjoint action, i.e. for any $a\in A$ we have
	$$a\cdot  (S\o id)(W)=\v_t(a)(S\o id)(W).$$
\end{proposition}
\begin{proof}
	Denote $W=r_2R_1\o r_1R_2$, here the letter $r$ is a copy of $R$. We compute
	\begin{align*}
	a\cdot (S\o id)(W)&= h_{(1)}S(r_2R_1)S(h_2)\o h_{(3)}r_1R_2S(h_{(4)})\\
	 &= h_{(1)}S(R_1)S(r_2h_{(3)})\o r_1h_{(2)}R_2S(h_{(4)})\\
	 &= \v_t(h_{(1)})S(r_2R_1)\o r_1R_2\v_t(h_{(2)})\\
	 &= \v_t(h)S(r_2R_1)\o r_1R_2.
	\end{align*}
	In the second and third equalities we use condition (3) of Definition \ref{QT}. In the last equality we use formula (\ref{r3}).
\end{proof}

\subsection{The quasitriangular structure on the Drinfeld double }

\begin{definition}
	Let $(A, B, \left\langle ,\right\rangle )$ be a WMHA pairing. Given $R\in M(A\o B)$,  $R$ is called canonical if
	$$\left\langle R , a\o b \right\rangle=\left\langle a, b\right\rangle$$
	for all $a\in A, b\in B$.
\end{definition}

Here we consider the extension from $\left\langle A\o B , A\o B \right\rangle$ to $\left\langle M(A\o B) , A\o B \right\rangle$. Such extensions will be used frequently in the following.

\begin{example}
	Let $G$ be any groupoid. Consider the WMHA pairing $\left\langle A, B \right\rangle $  in  Example \ref{EX1}. For any $g, h\in G$, let $\delta_p$ be the function that is 1 in $p$ and 0 else. Now the canonical element $R\in M(A\o B)$ is given by $R=\sum_{p\in G} \delta_p\o \lambda_p$.
	
\end{example}

\begin{proposition}
	Let $(A, B, \left\langle ,\right\rangle )$ be a WMHA pairing and $R\in M(A\o B)$, the following  are equivalent
	\begin{itemize}
		\item[(1) ] $R$ is canonical,
		\item[(2) ]  $(id\o \left\langle a, \cdot \right\rangle)R=a$ for any $a\in A$,
		\item[(3) ] $(\left\langle \cdot, b\right\rangle \o id)R=b$ for any $b\in B$.
	\end{itemize}
\end{proposition}
\begin{proof}
Let $R=R_1\o R_2$, for any $a'\in A, b'\in B$,
$$\left\langle   (\left\langle \cdot, b\right\rangle \o id)R , b'\right\rangle=\left\langle \left\langle R_1, b  \right\rangle R_2, b' \right\rangle=\left\langle R, b\o b' \right\rangle.$$
It implies that $(1)\Leftrightarrow (3)$. Note that if we replace $b$ by $a\rhd b$ and $b'$ by $a'\rhd b'$, then everything is covered. The proof of equivalence between $(1)$ and $(2)$ is similar.
\end{proof}

\begin{remark}
Given a WMHA pairing $(A, B, \left\langle ,\right\rangle )$, we can construct the Drinfeld double $D=A\o_D B$, and we have  non-degenerate algebra embeddings
$$A\hookrightarrow M(D),  \quad    B\hookrightarrow M(D)$$
  The two embeddings  give rise to the following
  $$A\o A \hookrightarrow M(A\o D),\quad  B\o B\hookrightarrow M(D\o B)$$
  $$A\o B\hookrightarrow M(D\o B),\quad  A\o B\hookrightarrow M(A\o D).$$
 Moreover, we have
 $$A\o A\hookrightarrow M(D\o D), a\o a' \mapsto (a\o_D 1)\o (a'\o_D  1)$$
 $$B\o B\hookrightarrow M(D\o D), b\o b'\mapsto (1\o_D b)\o (1\o_D b').$$
 Similarly, $A\o B$ can be viewed as a subalgebra of $M(D\o D)$, now the embedding is $a\o b \mapsto (a\o_D 1)\o (1\o_D b)$. Moreover they can be extended to the multiplier algebra. For any  $R\in M(A\o B)$, the image of $R$ in $M(D\o D)$ will be denoted by the same symbol $R$ in the sequel.
\end{remark}

\begin{proposition}\label{eee}
	Let $R$ be a canonical element for the WMHA pairing $\left\langle A, B \right\rangle $, then we have
	\begin{itemize}
		\item [(1)]  $(\D_A\o id_B)R=R^{13}R^{23} \quad \text{and} \quad (id_A\o \D_B)R=R^{12}R^{13}$
		\item[(2)] $(id_A\o \v_B)R=1 \quad \text{and} \quad (\v_A\o id_B)R=1$.
	\end{itemize}
\end{proposition}
\begin{proof}
First note that these equations are meaningful in the multiplier algebras.

(1) Denote $N=(\D_A\o id_B)(R), \quad Q=R^{13}R^{23}$ and $f= \left\langle \cdot, b\o b' \right\rangle\o id_B$, here $N,Q\in M(A\o A\o B)$. For any $b, b'\in B$, we have
$$f(N)=(\left\langle \cdot, bb'\right\rangle\o id_B)R=bb'$$
$$f(Q)=((\left\langle \cdot, b\right\rangle\o id_B)R)((\left\langle \cdot, b'\right\rangle\o id_B)R)=bb',$$
then we can obtain $N=Q$. The proof of second equation  is similar.

(2) Applying $id_A\o id_B\o \v_B$ on two sides of  $(id_A\o \D_B)R=R^{12}R^{13}$, then we can get $(id_A\o \v_B)R=1$. Similarly, we can prove the second formula.
\end{proof}

\begin{lemma}\label{ccc}
	Let $D$ be the Drinfeld double associated with WHMA pairing $\left\langle A, B \right\rangle$. For any $a, x\in A, b, y\in B$, we have
	$$\sum \left\langle a_{(2)}, b_{(1)}\right\rangle (1\o_D b_{(1)})(a_{(2)}x\o_D y)=\sum (1\o_D b_{(1)})(a_{(2)}x\o_D y)\left\langle a_{(1)}, b_{(2)}\right\rangle.$$
\end{lemma}
\begin{proof}
	\begin{align*}
	RHS &= a_{(3)}x_{(2)}\o_D b_{(2)}y \left\langle a_{(4)}x_{(3)}, b_{(1)}  \right\rangle \left\langle a_{(2)}x_{(1)}, S^{-1}(b_{(3)})  \right\rangle \left\langle a_{(1)}, b_{(4)}  \right\rangle\\
	&= a_{(3)}x_{(2)}\o_D b_{(3)}y\left\langle a_{(4)}, b_{(1)} \right\rangle \left\langle x_{(3)}, b_{(2)} \right\rangle \left\langle a_{(2)}, S^{-1}(b_{(5)}) \right\rangle \left\langle x_{(1)}, S^{-1}(b_{(4)}) \right\rangle \left\langle a_{(1)}, b_{(6)} \right\rangle\\
	&= (a_{(3)}\o_D b_{(2)})(x\o_D y) \left\langle a_{(4)}, b_{(1)} \right\rangle \left\langle a_{(2)}, S^{-1}(b_{(3)}) \right\rangle \left\langle a_{(1)}, b_{(4)} \right\rangle\\
	&= (a_{(1)}\lhd \v_s'(b_{(3)})\o_D b_{(2)})(x\o_D y)\left\langle a_{(2)}, b_{(1)} \right\rangle\\
	&\stackrel{(\ref{4}}{=} \left\langle a_{(1)}, b_{(2)}\right\rangle (1\o_D b_{(1)})(a_{(2)}x\o_D y)\\
	&=LHS.
	\end{align*}
\end{proof}

For convenience, only in this subsection the triple data $(\left\langle A, B \right\rangle, R)$ means that $R$ is a canonical element for the WMHA pairing $\left\langle A, B \right\rangle $.

\begin{proposition}\label{fff}
	Given  $(\left\langle A, B \right\rangle, R)$. For any $a\in A, b\in B$, we have
	\begin{itemize}
		\item [(1)] $R\D(a)=\D^{cop}(a)R$
		\item[(2)] $R\D^{cop}(b)=\D(b)R$
	\end{itemize}
\end{proposition}
\begin{proof}
	For any $a, x\in A, b, y\in B$,
	\begin{align*}
	(\left\langle \cdot, b \right\rangle\o id_D)(R\D(a)(1\o (x\o_D y)))
	&=\left\langle R_1a_{(1)}, b\right\rangle (1\o_D R_2)(a_{(2)}x\o_D y)\\
	&= \left\langle a_{(1)}, b_{(2)}\right\rangle (1\o_D b_{(1)})(a_{(2)}x\o_D y),
	\end{align*}
	On the other side, we have
	\begin{align*}
	(\left\langle \cdot, b \right\rangle\o id_D)(\D^{cop}(a)R(1\o (x\o_D y)))&= \left\langle a_{(2)}R_1, b\right\rangle (a_{(1)}\o_D R_2)(x\o_D y)\\
	&= (a_{(1)}\o_D b\lhd a_{(2)})(x\o_D y)
	\end{align*}
	By Lemma \ref{ccc}, we get
	$$(\left\langle \cdot, b \right\rangle)(R\D(a)(1\o (x\o_D y)))= (\left\langle \cdot, b \right\rangle\o id_D)(\D^{cop}(a)R(1\o (x\o_D y))).$$
	Since the pairing and product are non-degenerate, we have $R\D(a)=\D^{cop}(a)R$. The second equation can be proved in a similar way. We have unital action here, those elements can be covered.
\end{proof}

Now we can prove that $R$ satisfies the Yang-Baxter  equation in $M(A\o D\o B)$, here $D=A\o_D B$ is the Drinfeld double.

\begin{proposition}
	Given  $(\left\langle A, B \right\rangle, R)$, $R$ satisfies
	$$R^{23}R^{13}R^{12}=R^{12}R^{13}R^{23}.$$
\end{proposition}
\begin{proof}
	By Proposition \ref{eee} and \ref{fff}, we have
	\begin{align*}
	R^{12}R^{13}R^{23} &= (R\o 1)(\D \o id)(R)\\
	&= (\D^{cop}_A\o id)(R)R^{12}\\
	&= R^{23}R^{13}R^{12}.
	\end{align*}
\end{proof}

Finally we will show that $R$ is a quasitriangular structure in the sense of Definition \ref{QT}.

\begin{proposition}\label{E-inverse}
 Let	$\left\langle A, B \right\rangle$ be a WMHA pairing. If $R\in M(A\o B)$ is the canonical element and $D=A\o_D B$ is the Drinfeld double, then the image of $R$ is $E$-inverse in $M(D\o D)$.
\end{proposition}
\begin{proof}
	Set $R=u\o v$, the image of $R$ is $u\o_D 1\o 1\o_D v$ and we still denote by $R$. Define $\bar{R}=S(u)\o v$. By definition we have to show that $\bar{R}R=E_D, R\bar{R}=E_D^{cop}$ and $RE_D=R, \bar{R}E^{cop}=\bar{R}$. The canonical idempotents in $A$ and $B$ are denoted by $E$ and $e$ respectively.
	We compute
	\begin{align*}
	\bar{R}R &= (S(u¡¯)\o_D 1\o 1\o_D v')(u\o_D 1\o 1\o_D v)\\
	  &=S(u')u\o_D 1\o 1\o_D v'v\\
	  &=\v_s(u)\o_D 1\o 1\o v
	\end{align*}
	Now $\bar{R}R=E_D$ is equivalent to
	$$\v_s(u)\o_D 1\o 1\o v=S^{-1}(e_2)\rhd E_1\o_D 1 \o 1 \o_D E_2\rhd e_1,$$
	which can be regarded as an equality in $B\o A$:
	$$\v_s(u)\o v=\left\langle E_2, S^{-1}(e_3) \right\rangle E_1\o e_1\left\langle E_3, e_2\right\rangle.$$
	Evaluating both sides on any $a\in A$, we have
	\begin{align*}
	\v_s(a)&= \v_A(aE_3S^{-1}(E_2))E_1=\v_A(aE_2)E_1.
	\end{align*}
	The proof of $R\bar{R}=E_D^{cop}$ is similar.  Next let us compute  $RE_D=R$.
	\begin{align*}
	RE_D &=(u\o_D 1\o 1\o_D v)(E_1\o_D e_2\o E_2\o_D e_1)\\
	        &=u(S^{-1}(e_2)\rhd E_1)\o_D 1\o 1\o_D v(E_2\rhd e_1)
	\end{align*}
	So it is equivalent to show that
	$$u(S^{-1}(e_2)\rhd E_1)\o v(E_2\rhd e_1)=u\o v.$$
	Taking arbitrary $b\in B$, we get
	\begin{align*}
	b &=\left\langle u(S^{-1}(e_2)\rhd E_1), b \right\rangle v(E_2\rhd e_1) \\
	&= \v_B(b_{(2)}S^{-1}(e_3)e_2)e_1\\
	&=\v_B(b_{(2)})b_{(1)}.
	\end{align*}
	The second formula is similar.
\end{proof}

\begin{theorem}\label{QT-main}
  Given a weak multiplier Hopf algebra pairing $\left\langle A, B \right\rangle$. Let $D=A\o_D B$ be the Drinfeld double and $R$ be the canonical element. Then $(D, R)$ is a quasitriangular weak multiplier Hopf algebra.
\end{theorem}
\begin{proof}
	Note that $R\in M(A\o B)$ can be embedded into $M(D\o D)$.  Let $R=R_1\o R_2=r_1\o r_2, E=E_1\o E_2=e_1\o e_2$, we have
	\begin{align*}
	(\D_D\o id)(R)&=(\D_D\o id)((R_1\o_D 1)\o (1\o_D R_2))\\
	&= R_{1(1)}\o_D E_2 \o R_{1(2)}\o_D E_1\o 1\o_D R_2\\
	&= R_{1(1)}\o_D 1 \o R_{1(2)}\o_D 1\o  1\o_D R^2\\
	&= R_1\o_D 1\o r_1\o_D 1\o 1\o_D R_2r_2\\
	&= R^{13}R^{23}
	\end{align*}
	In the third equality we use $RE=R$.  For the fourth equality we use $(\D\o id)R=R^{13}R^{23}$.  Similarly, we obtain that
	$(id\o \D_D)(R)=R^{13}R^{12}$.
	
	Finally, we have
	$$\D^{cop}(a\o_D b)R=\D^{cop}(a)\D(b)R=\D^{cop}(a)R\D^{cop}(B)=R\D(a\o_D b).$$
	With Proposition \ref{E-inverse}, we finish the proof.
\end{proof}

Generally it is not easy to  find a canonical element. But we can find one in the following special case. Consider an algebraic quantum groupoid $A$ and its dual $\widehat{A}$ with the natural pairing. Define multiplier $R$ in $M(\widehat{A}\o A)$ by
$$
R(\varphi(a\cdot)\o a')=\sum \varphi(a_{(2)}\cdot)\o S(a_{(1)})a'
$$
$$
(\psi(\cdot a)\o a')R=\sum \psi(\cdot a_{(1)})\o a'S(a_{(2)}),
$$
here $a, a'\in A$, $\varphi$ and $\psi$ are left and right  integrals on $A$, respectively.

\begin{proposition}
	Take the notations as above. $R$ is canonical for $\left\langle \widehat{A}, A \right\rangle$.
\end{proposition}
\begin{proof}
	For any $b\in A, f\in \widehat{A}$, we can choose $\varphi(a\cdot)\in \widehat{A}$ and $a'\in A$ such that $\varphi(a\cdot)\rhd b=b$ and $a'\rhd f=f$. Now we compute
	\begin{align*}
	\left\langle R, f\o b \right\rangle
	&= \left\langle R(\varphi(a\cdot)\o a'), f\o b \right\rangle\\
	&= \left\langle \varphi(a_{(2)}\cdot)\o S(a_{(1)})a', f\o b \right\rangle\\
	&=\left\langle f,\varphi(a_{(2)}b)S(a_{(1)})a'\right\rangle\\
	&=\left\langle f,\varphi(ab_{(2)})b_{(1)}a'\right\rangle\\
	&= \left\langle f, b \right\rangle.
	\end{align*}
	In the fourth equality we use Proposition 1.5 of \cite{VW4} which states that $\varphi$ is a left integral if and only if
	$$(id\o \varphi)((1\o a)\D(b))=S((id\o \varphi)(\D(a)(1\o b)))$$
	for all $a$ and $b$.
\end{proof}

Under certain conditions, we can give the concrete expression of $R$. First we need the concept of cointegral in a weak multiplier Hopf algebra.

\begin{definition}
	Let $A$ be a  weak multiplier Hopf algebra.
	A non-zero element  $h\in M(A)$ is called a {\it left cointegral} if it satisfies
	$ah=\v_t(a)h $
	for all $a\in A$. Similarly, a  non-zero element  $h\in M(A)$ is called a {\it right cointegral} if it satisfies
	$ka=k\v_s(a)$
	for all $k\in A$.
\end{definition}

And we also have that $ah=\v_t(a)h$ is equivalent with $(1\o a)\D(h)=(S(a)\o 1)\D(h)$. The details about cointegrals will be discussed in another paper.

\begin{proposition}
	Let $A$ be an algebraic quantum groupoid with a unit. If $\varphi$ is a left integral on $A$ and $t$ a left cointegral in $A$ such that $\varphi(t)=1$, then the element
	$$\sum \varphi(\cdot t_{(2)})\o S^{-1}(t_{(1)})\in \widehat{A}\o A$$
	is canonical for $\left\langle \widehat{A}, A \right\rangle$.
\end{proposition}
\begin{proof}
	Again we use Proposition 1.5 of \cite{VW4}, for any $f\in \widehat{A}, a\in A$, we have
	\begin{align*}
	\left\langle  \varphi(\cdot t_{(2)})\o S^{-1}(t_{(1)}), f\o a \right\rangle
	&= \left\langle \varphi(at_{(2)})f, S^{-1}(t_{1}) \right\rangle\\
	&= \left\langle f, \varphi(t_{(2)})S^{-1}(S(a)t_{(1)}) \right\rangle\\
	&=\left\langle f, a \right\rangle.
	\end{align*}
\end{proof}

 {\bf Acknowledgment}

 The work of Nan Zhou was  supported by  Zhejiang Shuren University grant KXJ1418601 and KBY0119604C/014.


\begin{thebibliography}{100}
\bibitem{AT} Adi Ben-Israel, Thomas N.E. Greville, Generalized Inverses Theory and Applications, second edition. CMS books in mathematics(15).
\bibitem{B1} G. B\"{o}hm, Doi-Hopf modules over weak Hopf algebras, Comm. Algebra 28(10) (2010) 4687--4698.
\bibitem{B2} G. B\"{o}hm, Comodules over weak multiplier bialgebras, Int. J. Math. 25(5) (2014) 1450037.
\bibitem{B3} G. B\"{o}hm, Yetter-Drinfeld modules over weak multiplier bialgebras, Israel J. Math. 209(1) (2015) 85--123.
\bibitem{BGL} G. B\"{o}hm, J. Gomez-Torrecillas, E. Lopez-Centella, Weak multiplier bialgebras, Trans. Amer. Math. Soc. 367(12) (2015) 8681--8721.
\bibitem{BNS}  G. B\"ohm, F. Nill and K. Szlach\'{a}nyi, Weak Hopf algebras I. Integral theory and $C^*$-structure,  J. Algebra  221 (1999) 385--438.
\bibitem{BS}  G. B\"{o}hm, K. Szlach\'{a}nyi, Weak Hopf algebras II. Representation theory and
 the Markov trace,  J. Algebra 233 (2000) 156--212.
\bibitem{De1} L. Delvaux, Twisted tensor product of multiplier Hopf (*-)algebras, J. Algebra 269 (2003) 285--316.
\bibitem{DeVD}  L. Delvaux, A. Van Daele, The Drinfeld double of multiplier Hopf algebras, J. Algebra 272(1) (2004) 273--291.
\bibitem{DeVD1} L. Delvaux, A. Van Daele, The Drinfeld double versus the Heisenberg double for an algebraic quantum group, Journal of Pure and Applied Algebra 190 (2004) 59--84.
\bibitem{DeVDW}L. Delvaux, A. Van Daele, S. Wang, Quasitriangular (G-cograded) multiplier Hopf algebras, J. algebra 289(2005) 484--514.
\bibitem{DrVD} B. Drabant, A. Van Daele, Pairing and Quantum Double of Multiplier Hopf Algebras, Algebra and Representation Theory  4 (2001) 109--132.
\bibitem{Lu} J.-H. Lu, On the Drinfeld double and the Heisenberg double of a Hopf algebra, Duke Math. J. 74(1994) 763-776.
\bibitem{Ma} S. Majid, Foundations of Quantum Group Theory, Cambridge University Press, Cambridge, 1995.
\bibitem{M} S. Montgomery, Hopf algebras and their actions on rings, AMS, Providence, 1993.
\bibitem{NTV} D. Nikshych, V. Turaev, L. Vainerman, Invariants of knots and 3-manifolds from quantum groupoids, Topology and its Applications 127 (2003) 91--123.
\bibitem{VD} A. Van Daele, Multiplier Hopf algebras, Trans. Amer. Math. Soc. 342(2) (1994) 917--932.
\bibitem{VD1} A. Van Daele, An algebraic framework for group duality,  Adv. Math. 140 (1998) 323--366.
\bibitem{VD2} A. Van Daele, Tools for working with multiplier Hopf algebras. Arab. J. Sci. Eng. 33(2C), (2008) 505-527.
\bibitem{VW1} A. Van Daele and S. Wang, Weak Multiplier Hopf algebras. Preliminaries, motivation and basic examples,  Operator Algebras and Quantum Groups. W.Pusz and P.M. Soltan (eds), Banach Center Publications (Warsaw), 98 (2012) 367--415.
\bibitem{VW2} A. Van Daele and S. Wang,Weak Multiplier Hopf Algebras. The main theory, Journal f\"{u}r die reine und angewandte Mathematik (Crelle's Journal), 705 (2015) 155--209.
\bibitem{VW3} A. Van Daele and S. Wang, Weak Multiplier Hopf Algebras II. The source and target algebras, Symmetry 12 (12) (2020),1975
\bibitem{VW4} A. Van Daele, S.  Wang, Weak multiplier Hopf algebras III. Integrals and Duality, arXiv: 1701.0495v3.
\bibitem{WZW} W. Wang, N. Zhou, S. Wang, Semidirect products of weak multiplier Hopf algebras: Smash products and smash coproducts.  Comm. Algebra, 46(2018) 3241--3261.
\bibitem{Zhang} Y. Zhang, The quantum double of a coFrobenius Hopf algebra, Comm. Algebra, 27:3, (1999) 1413-1427.
\bibitem{ZhuW} H. Zhu, S. Wang, A generalized Drinfeld quantum double construction based on weak Hopf algebras, Comm. Algebra 38(1) (2009) 199--229.
\bibitem{ZhouW} N. Zhou, S. Wang, A duality theorem for weak multiplier Hopf algebras actions. Int. J. Math. 28(5) (2017) 1750032.
\end{thebibliography}
\end{document}